\documentclass[preprint,12pt]{elsarticle}

\usepackage{amssymb}

\usepackage{amsmath, amssymb,latexsym}
\usepackage{amsthm}
\usepackage{graphicx}

\theoremstyle{definition}
\newtheorem{thm}{Theorem}[section]
\newtheorem*{thm*}{Theorem}
\newtheorem{defi}[thm]{Definition}
\newtheorem{lemm}[thm]{Lemma}

\newtheorem{remark}[thm]{Remark}
\newtheorem{coro}[thm]{Corollary}

\journal{Differential Geometry and its Applications}

\begin{document}

\begin{frontmatter}



\title{On cohomogeneity one hyperpolar actions related to $G_{2}$}


\author[1]{Shinji Ohno}

\ead{ohno.shinji@nihon-u.ac.jp}

\affiliation[1]{organization={Department of Mathematics College of Humanities and Sciences, Nihon University},
            addressline={3-25-40, Sakurajosui}, 
            city={Setagaya-ku},
            postcode={156-8550}, 
            state={Tokyo},
            country={Japan}}

\author[2]{Yuuki Sasaki}

\ead{y_sasaki@cc.utsunomiya-u.ac.jp}

\affiliation[2]{organization={Cooperative Faculty of Education, Utsunomiya University},
            addressline={350, Mine-Machi}, 
            city={Utsunomiya},
            postcode={321-8505}, 
            state={Tochigi},
            country={Japan}}

\begin{abstract}

Cohomogeneity one actions on irreducible Riemannian symmetric spaces of compact type are classified into three cases: Hermann actions, actions induced by the linear isotropy representation of a Riemannian symmetric space of rank 2, and exceptional actions.
In this paper, we consider exceptional actions related to the exceptional compact Lie group $G_{2}$ and investigate some properties of their orbits as Riemannian submanifolds.
In particular, we examine the principal curvatures of principal orbits and classify principal orbits that are minimal, austere, weakly reflective, and proper biharmonic.

\end{abstract}




\end{frontmatter}



\section{Introduction}

The theory of constructing minimal submanifolds using Lie group actions was introduced by Hsiang and Lawson and they classified homogeneous minimal hypersurfaces in the unit sphere (\cite{HsL}).
Such hypersurfaces are some orbits of a linear isotropy action of a Riemannian symmetric space of rank 2.
Since then, many geometers have developed this theory and studied orbits of various isometry group actions of Lie groups as Riemannian submanifolds. 
For example, the orbits of a linear isotropy action of a Riemannian symmetric space are called {\it $R$-spaces}, and $R$-spaces are naturally embedded in round spheres.
The properties of $R$-spaces as Riemannian submanifolds of round spheres have been studied extensively (\cite{HTST}, \cite{KO}, etc.).
An isometric action of a compact Lie group on a Riemannian manifold is called a {\it polar} action if there exists a closed and connected submanifold $\Sigma$ that intersects all orbits orthogonally. 
Such $\Sigma$ is called a {\it section}.
A section is automatically a totally geodesic submanifold.
If the section is also flat, then the action is called a {\it hyperpolar} action.
Let $G$ be a compact Lie group and $\theta_{i}\ (i=1,2)$ be involutive automorphisms of $G$.
Set $F(\theta_{i}, G) = \{ g \in G \ ;\ \theta_{i}(g) = g \}$ and denote the identity component of $F(\theta_{i}, G)$ by $F_{0}(\theta_{i}, G)$.
Let $K_{i}$ be closed subgroups of $G$ such that $F_{0}(\theta_{i}, G) \subset K_{i} \subset F(\theta_{i}, G)$.
Then, $G/K_{i}$ are compact Riemannian symmetric spaces and the actions $K_{i} \times G/K_{j} \rightarrow G/K_{j} ; (k_{j}, gK_{i}) \mapsto k_{j}gK_{i}\ (i, j = 1,2)$ are called {\it Hermann actions}.
It is well known that any Hermann action is a hyperpolar action.
Moreover, any cohomogeneity one action on an irreducible Riemannian symmetric space of compact type is hyperpolar.
Kollross classified cohomogeneity one actions on irreducible Riemannian symmetric spaces of compact type up to the orbit equivalence and such actions can be categorized into three cases (\cite{K1}), namely

\begin{itemize}

\item[(i)]
cohomogeneity one Hermann actions,

\vspace{-3mm}

\item[(ii)]
actions induced by the linear isotropy representation of a Riemannian symmetric space of rank 2,

\vspace{-3mm}

\item[(iii)]
exceptional actions.

\end{itemize}

\hspace{-22pt}
To classify homogeneous minimal hypersurfaces of irreducible Riemannian symmetric spaces of compact type, we need to study principal orbits of the above actions.
For Hermann actions, a method to compute the second fundamental form of orbits was given by Ikawa in the cases of $\theta_{1}\theta_{2} = \theta_{2}\theta_{1}$ (\cite{I1}) and by the first author in the case of $\theta_{1}\theta_{2} \not= \theta_{2}\theta_{1}$ (\cite{Ohno}).

In this paper, we consider the third cases, namely, exceptional actions.
According to the result of Kollross (\cite{K1}), there are seven distinct types.
Table \ref{tab:my_label} lists all exceptional actions.
In this table, $H$ and $K$ are subgroups of $G$ and the following actions are cohomogeneity one hyperpolar actions:
\[
\begin{split}
& (H \times K) \times G \rightarrow G \ ;\ ( (h,k), g) \mapsto hgk^{-1}, \\
& H \times G/K \rightarrow G/K\ ;\ (h, gK) \mapsto hgK, \\
& K \times G/H \rightarrow G/H \ ;\ (k, gH) \mapsto kgH.
\end{split}
\]
In \cite{Verhoczki}, Verhoczki studied the $H \times K$-action on $G$ of type (I) and the $K$-action on $G/H$ of type (II).
Moreover, Enoyoshi studied the $K$-action on $G/H$ of type (IV) \cite{E}.
In this paper, we focus on the remaining cases related to the exceptional compact Lie group $G_{2}$, namely, types (I), (II), (III), (IV), and (V).
Note that in these types, both $G/K$ and $G/H$ are compact symmetric spaces.
In particular, $G/K = S^{6}$ in types (I) and (II), and $G/K = \mathbb{R}P^{7}$ in types (III), (IV), and (V).
Moreover, $SO(7)/U(3)$ is isomorphic to some oriented real Grassmannian.
The $H$-action on $G/K$ of types (I), (II), (III), (IV), and (V) and the $K$-action on $G/H$ of types (I), (III), and (V) are orbit equivalent to Hermann actions.
Therefore, we consider the $H \times K$-action on $G$ of types (II), (III), (IV), and (V).

\begin{table}[h]
    \centering
    \begin{tabular}{|c|c|c|c|c|c|}\hline
     type    & $G$ & $H$   & $K$  \\ \hline \hline
     (I)   & $G_{2}$          & $SU(3)$                         &  $SU(3)$                   \\ \hline
     (II)  & $G_{2}$          & $SO(4)$                         &  $SU(3)$                   \\ \hline
     (III)  & $SO(7)$ &$G_{2}$                                 &  $G_{2}$ \\ \hline
     (IV)   & $SO(7)$ &$SO(4)\times SO(3)$  &   $G_{2}$ \\ \hline
     (V)  &$SO(7)$ &$U(3)$                          &$G_{2}$  \\ \hline
     (VI)  & $SO(16)$&$SO(14)\times SO(2)$   &  $Spin(9)$ \\ \hline
     (VII)  & $SO(4n)$&$SO(4n-2)\times SO(2)$   &  $Sp(n)Sp(1)$ \\ \hline
    \end{tabular}
    \caption{Exceptional cohomogeneity one actions}
    \label{tab:my_label}
\end{table}

This paper is organized as follows:
In Section 2, we provide preliminary results.
We recall some fundamental properties related to group actions, submanifold theory, and octonions.
Additionally, we consider some subgroups of $G_{2}$ and $SO(7)$.
In Section 3, we analyze type (II) and calculate the principal curvatures of principal orbits.
We show that there exists exactly one minimal principal orbit, which is neither austere nor weakly reflective.
Furthermore, we demonstrate that there are exactly two proper biharmonic principal orbits.
Section 4 focuses on type (III), while Section 5 addresses type (IV).
In both cases, we find that there exists exactly one minimal principal orbit, which is weakly reflective.
We compute the principal curvatures of principal orbits and show that there exists exactly one proper biharmonic principal orbit in type (III) and exactly two proper biharmonic principal orbits in type (IV).
In Section 6, we examine type (V).
Here, we show that there exists exactly one minimal principal orbit, which is neither austere nor weakly reflective.
We compute the principal curvatures of principal orbits and demonstrate that there exist exactly two proper biharmonic principal orbits.
Moreover, in type (IV), each $K$-orbit in $G/H$ is characterized as a level set of a $K$-invariant function on $G/H$ (\cite{HL}, \cite{E}).
By using this function, both $K \times H$-orbits in $G$ and $H$-orbits in $G/K$ can be characterized in a similar way.
For types (II), (III), and (V), we also consider such characterization of $K$-orbits in $G/H$.


\section{Preliminaries}

\subsection{Some submanifolds}

Let $(N,h)$ be a complete Riemannian manifold and $M$ be a submanifold of $N$.
For any $x \in M$, we denote by $T^{\perp}_{x}M$ the normal subspace of $T_{x}M$ in $T_{x}N$.
Furthermore, for any $\xi \in T^{\perp}_{x}M$, we denote the shape operator by $A^{\xi}$.

\begin{defi} \cite{HL}
$M$ is an {\it austere} submanifold of $N$ if for any $x \in M$ and $\xi \in T_{x}^{\perp}M$, the set of eigenvalues of $A^{\xi}$ with multiplicities is invariant under multiplication by $-1$.
\end{defi}

By definition, any austere submanifold is also a minimal submanifold.
Harvey and Lawson demonstrated that a special Lagrangian cone in a complex Euclidean space can be constructed using an austere submanifold of the unit sphere \cite{HL}.
Ikawa, Sakai, and Tasaki introduced a special class of austere submanifolds in \cite{IST2}.

\begin{defi} \cite{IST2}
$M$ is a {\it weakly reflective} submanifold of $N$ if for any $x \in M$ and $\xi \in T^{\perp}_{x}M$, there exists an isometry $\sigma_{\xi}$ of $N$ such that $\sigma_{\xi}(x) = x, d\sigma_{\xi}(\xi) = -\xi$, and $\sigma_{\xi}(M) \subset M$.
\end{defi}

If $M$ is a connected component of the fixed point set of an involutive isometry of $N$, then we call $M$ a {\it reflective} submanifold.
It is clear that a reflective submanifold is a weakly reflective submanifold.
However, the converse does not hold.
While any reflective submanifold is totally geodesic, a weakly reflective submanifold is not necessarily totally geodesic.
Ikawa, Sakai, and Tasaki proved that every weakly reflective submanifold is an austere submanifold and established the following fundamental properties of weakly reflective submanifolds.

\begin{lemm}\cite{IST2} \label{wr1}
Let $G$ be a connected Lie group acting isometrically on $(N, h)$.
Denote the $G$-orbit through $x \in N$ by $G(x)$.
If for each $\xi \in T_{x}^{\perp}G(x)$, there exists an isometry $\sigma_{\xi}$ of $N$ such that $\sigma_{\xi}(x) = x, d\sigma_{\xi}(\xi) = -\xi,$ and $\sigma_{\xi}(G(x)) = G(x)$, then $G(x)$ is a weakly reflective submanifold.
\end{lemm}

\begin{lemm}\cite{IST2} \label{wr2}
Any singular orbit of a cohomogeneity one action on a Riemannian manifold is a weakly reflective submanifold.
\end{lemm}

Next, we recall harmonic maps and biharmonic maps.
Let $(M,g)$ be a compact Riemannian manifold and $(N,h)$ be a Riemannian manifold.
Denote by $\nabla^{M}$ and $\nabla^{N}$ the Levi-Civita connection of $M$ and $N$, respectively. 
A smooth map $\phi : M \rightarrow N$ is a {\it harmonic map} if $\phi$ is an extremal of the energy functional defined by
\[
E(\phi) = \frac{1}{2} \int_{M} |d\phi|^{2} v_{g}.
\]
Let $\phi^{*}TN$ be the pullback bundle of the tangent bundle $TN$ induced by $\phi$.
For any vector field $V$ on $N$, $(V \circ \phi)(x) = V_{\phi(x)}$ defines a section of $\phi^{*}TN$.
There exists a unique connection $\nabla^{\phi}$ on $\phi^{*}TN$ such that 
\[
\Big( \nabla^{\phi}_{X} (V \circ \phi) \Big)_{x} = \Big( \nabla^{N}_{d\phi(X)}V \Big)_{\phi(x)},
\]
where $x \in M$, $X$ is a vector field on $M$, and $V$ is a vector field on $N$.
We define the second fundamental form $\nabla(d\phi)$ of $\phi$ as follows:
For any $x \in M$ and vector fields $X$ and $Y$ on $M$,
\[
\nabla(d\phi)(X,Y) = \nabla^{\phi}_{X} d \phi(Y) - d\phi \Big( \nabla_{X}^{M}Y \Big).
\]
Then, $\tau(\phi) =  \mathrm{tr} \nabla(d\phi)$ is called the {\it tension field} of $\phi$.
It is well known that $\phi$ is a harmonic map if and only if $\tau(\phi) = 0$.
Moreover, $\phi$ is a {\it biharmonic map} (\cite{Eelles-Lemaire}) if $\phi$ is an extremal of the bienergy functional defined by
\[
E_{2}(\phi) = \frac{1}{2} \int_{M} |\tau(\phi)|^{2} v_{g}.
\]
By definition, any harmonic map is a biharmonic map.
If a biharmonic map is not harmonic, then it is called a {\it proper} biharmonic map.
In the following, we consider isometric embeddings.
Then, the tension field is the mean curvature vector field.

\begin{lemm} \cite{Ou} \label{biharmonic}
Let $(N,h)$ be a Riemannian manifold with $\dim N \geq 3$.
Denote the Ricci tensor of $N$ by $r$.
Suppose that $N$ is Einstein, i.e., there exists $\lambda \in \mathbb{R}$ such that $r = \lambda h$.
Let $M$ be a hypersurface of $N$.
Assume that the mean curvature vector field of $M$ is constant.
Then, $M$ is a biharmonic submanifold if and only if $M$ is either minimal or non-minimal with $|A^{\xi}|^{2} = \lambda$, where $\xi$ is a unit normal vector field of $M$.
\end{lemm}

Let $G$ be a semisimple compact Lie group and $H$ and $K$ be subgroups of $G$.
The $H \times K$-orbit through $x \in G$ is denoted by $HxK$.
In the following, the left and right translations by $g \in G$ are denoted by $L_{g}$ and $R_{g}$, respectively.
Then, 
\[
\begin{split}
T_{x}HxK &= \left\{ \left. \frac{d}{dt} \exp(tX) x \exp (-tY) \right|_{t=0} \  ;\  X \in \mathfrak{h}, Y\in \mathfrak{k} \right\} \\
&= \{ dL_{x} ( \mathrm{Ad}(x)^{-1}X - Y) \ ;\  X \in \frak{h}, Y \in \frak{k} \} \\
&= dL_{x} \big(\mathrm{Ad}(x)^{-1}\mathfrak{h} +\mathfrak{k} \big).
\end{split}
\]
Let $Z = (X,Y) \in \frak{g} \times \frak{g}$.
The Killing vector field on $G$ corresponding to $Z$ is denoted by $Z^{*}$.
That is, for any $x \in G$,
\[
Z^{*}_{x} = \left. \frac{d}{dt} \exp(tX) x \exp (-tY)\ \right|_{t=0} = dL_{x}(\mathrm{Ad}(x)^{-1} X -Y).
\]
We fix a bi-invariant Riemannian metric on $G$ and denote it by $\langle \ ,\ \rangle$.
Let $\nabla$ be the Levi-Civita connection of $(G,\langle \ ,\ \rangle)$.
By the Koszul formula, we obtain Lemma \ref{Levi-Civita}.

\begin{lemm} \label{Levi-Civita}
For any $x \in G$ and $Z_{1}=(X_{1}, Y_{1}), Z_{2}=(X_{2}, Y_{2})\in \mathfrak{g}\times \mathfrak{g}$, 
\[
\begin{split}
\left( \nabla_{Z_{1}^{\ast}}Z_{2}^{\ast}\right)_{x}
&= -\frac{1}{2} dL_{x} \Big[ \mathrm{Ad}(x)^{-1}X_{1} -Y_{1}, \mathrm{Ad}(x)^{-1}X_{2} +Y_{2} \Big] \\
&= -\frac{1}{2} dL_{x} \Big[ \mathrm{Ad}(x)^{-1}X_{1} -Y_{1}, \Big( \mathrm{Ad}(x)^{-1}X_{2} -Y_{2} \Big) + 2Y_{2} \Big]. \\
\end{split}
\]

\end{lemm}


\subsection{The compact Lie group $G_{2}$}

We recall several fundamental properties of octonions from \cite{Yokota}.
Let $\mathbb{O}= \sum_{i=0}^{7}\mathbb{R}e_{i}$ be the octonions, where $e_{0}, \cdots, e_{7}$ form a basis.
The multiplication between two octonions is defined as follows.
$e_{0}$ is the unit element and denoted by $1$.
For any $1 \leq i \not= j \leq 7$, we have $e_{i}^{2} = -1$ and $e_{i}e_{j} = -e_{j}e_{i}$.
Assume that this multiplication satisfies the distributive law.
In Figure $1$, the multiplication among $e_{1},e_{2},e_{3}$ is defined as $e_{1}e_{2} = e_{3}, \ e_{2}e_{3}=e_{1}, \ e_{3}e_{1}=e_{2}$.
Similarly, the multiplication among any three elements on each of the other lines and the circle is defined in the same manner. 
Note that the associative law does not follow in $\mathbb{O}$. 


\begin{figure}[htbp]
\begin{center}
\includegraphics[width = 50mm]{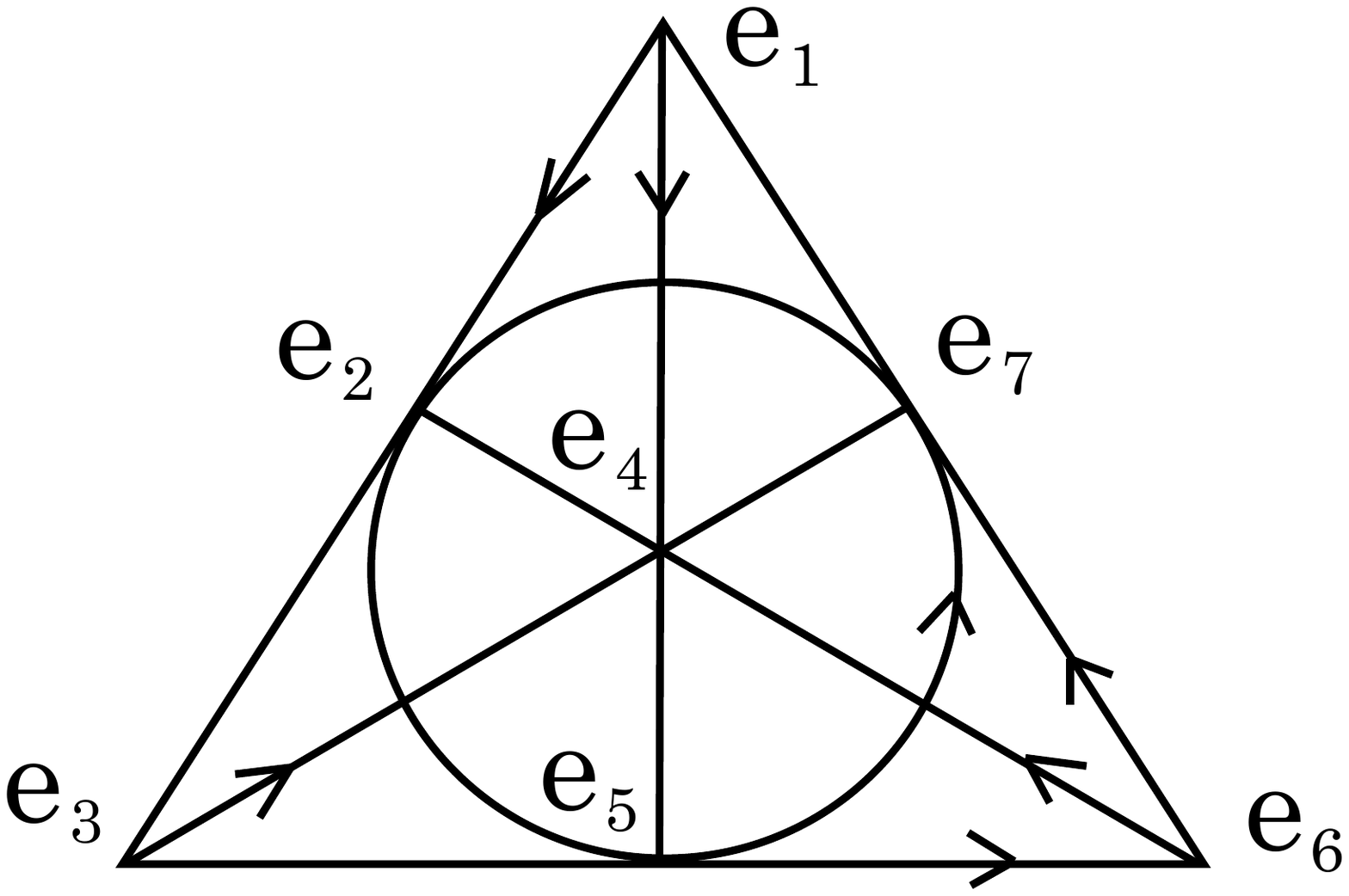}
\caption{Octonions}
\label{octonion}
\end{center}
\end{figure}

For any $x = x_{0} + \sum_{i=1}^{7}x_{i}e_{i} \in \mathbb{O} \ (x_{i} \in \mathbb{R})$, the element $\bar{x} = x_{0} - \sum_{i=1}^{7}x_{i}e_{i}$ is referred to as the conjugate of $x$.
Define $\mathrm{Im}\mathbb{O} = \{ x \in \mathbb{O} \ ;\ \bar{x} = -x \}$ and elements of $\mathrm{Im}\mathbb{O}$ are called pure octonions.
For any $x = \sum_{i=0}^{7}x_{i}e_{i}, y = \sum_{i=0}^{7}y_{i}e_{i} \in \mathbb{O}$, set $(x,y) = \sum_{i=0}^{7}x_{i}y_{i}$.
Then, $(\ ,\ )$ defines an inner product of $\mathbb{O}$, which satisfies
\[
\begin{split}
(x,y) &= (y,x) = \frac{1}{2}( \bar{x}y + \bar{y}x) = \frac{1}{2}(x\bar{y} + y\bar{x}).
\end{split}
\]
Set $|x| = \sqrt{(x,x)}$.
In the present paper, we denote $O(\mathbb{O})$ and $SO(\mathbb{O})$ by $O(8)$ and $SO(8)$, respectively.
Moreover, $O(7) = \{ g \in O(8)\ ;\ g(1) = 1 \}$ and $SO(7) = \{ g \in SO(8) \ ;\ g(1) = 1\}$.
If a linear automorphism $f:\mathbb{O} \rightarrow \mathbb{O}$ satisfies $f(xy) = f(x)f(y)$ for any $x,y \in \mathbb{O}$, then $f$ is called an automorphism of $\mathbb{O}$.
The group of all automorphisms of $\mathbb{O}$ is the exceptional compact Lie group $G_{2}$.
For any $g \in G_{2}$, it holds that $g(1) = 1$ and $(g(x), g(y)) = (x,y)$ for any $x,y \in \mathbb{O}$.
Moreover, $G_{2}$ is simply connected and $G_{2} \subset SO(7)$.
Let $\frak{so}(8)$ and $\frak{so}(7)$ denote the Lie algebras of $SO(8)$ and $SO(7)$, respectively.
Define the endomorphism $G_{ij} : \mathbb{O} \rightarrow \mathbb{O}$ for any $0 \leq i \not= j \leq 7$ by
\[
G_{ij} (e_{k} ) =
\left\{
\begin{array}{lllll}
e_{j} & (k = i), \\
-e_{i} & (k = j), \\
0 & (k \not= i,j ).
\end{array}
\right.
\]
It holds that $[G_{ik}, G_{kj}] = -G_{ij}$.
Moreover, $G_{ij}\ (0 \leq i < j \leq 7)$ form a basis of $\frak{so}(8)$, and $G_{ij}\ (1 \leq i < j \leq 7)$ form a basis of $\frak{so}(7)$.
Define the elements of $\frak{so}(7)$ as follows: 
\[
\begin{array}{lrlrlllll}
V_{1}(\lambda, \mu, \nu) = & \lambda G_{23} + \mu G_{45} + \nu G_{67}, & V_{2}(\lambda, \mu, \nu) = & -\lambda G_{13} - \mu G_{46} + \nu G_{57}, \\
V_{3}(\lambda, \mu, \nu) = & \lambda G_{12} + \mu G_{47} + \nu G_{56}, & V_{4}(\lambda, \mu, \nu) = & -\lambda G_{15} + \mu G_{26} - \nu G_{37}, \\
V_{5}(\lambda, \mu, \nu) = & \lambda G_{14} - \mu G_{27} - \nu G_{36}, & V_{6}(\lambda, \mu, \nu) = & -\lambda G_{17} - \mu G_{24} + \nu G_{35}, \\
V_{7}(\lambda, \mu, \nu) = & \lambda G_{16} + \mu G_{25} + \nu G_{34},
\end{array}
\]
where $\lambda, \mu, \nu \in \mathbb{R}$.
Set $V_{i} = \{ V_{i}(\lambda, \mu, \nu) \ ;\ \lambda, \mu, \nu \in \mathbb{R} \}$ for any $1 \leq i \leq 7$.
By direct computations, we obtain Lemma \ref{bracket}.

\begin{lemm} \label{bracket}
For any $1 \leq i, j, k \leq 7$, if $e_{i}e_{j} = \pm e_{k}$, then $[V_{i}, V_{j}] \subset V_{k}$.
In particular, 
\[
\begin{split}
& \big[ V_{1}(\lambda_{1}, \mu_{1}, \nu_{1}), V_{4}(\lambda_{4}, \mu_{4}, \nu_{4}) \big] = V_{5} \big( \mu_{1}\lambda_{4}, \ -(\lambda_{1}\nu_{4} + \nu_{1}\mu_{4}), \ -(\lambda_{1}\mu_{4} + \nu_{1}\nu_{4}) \big), \\
& \big[ V_{1}(\lambda_{1}, \mu_{1}, \nu_{1}), V_{5}(\lambda_{5}, \mu_{5}, \nu_{5}) \big] = V_{4} \big( -\mu_{1}\lambda_{5}, \ \lambda_{1}\nu_{5} + \nu_{1}\mu_{5}, \ \lambda_{1}\mu_{5} + \nu_{1}\nu_{5} \big), \\
& \big[ V_{4}(\lambda_{4}, \mu_{4}, \nu_{4}), V_{5}(\lambda_{5}, \mu_{5}, \nu_{5}) \big] = V_{1} \big( -(\mu_{4}\nu_{5} + \nu_{4}\mu_{5}), \ \lambda_{4}\lambda_{5}, \ -(\mu_{4}\mu_{5} + \nu_{4}\nu_{5}) \big), \\
& \big[ V_{2}(\lambda_{2}, \mu_{2}, \nu_{2}), V_{4}(\lambda_{4}, \mu_{4}, \nu_{4}) \big] = V_{6} \big( \lambda_{2}\nu_{4} + \nu_{2}\lambda_{4}, \ -\mu_{2}\mu_{4}, \ \lambda_{2}\lambda_{4} + \nu_{2}\nu_{4} \big), \\
& \big[ V_{2}(\lambda_{2}, \mu_{2}, \nu_{2}), V_{6}(\lambda_{6}, \mu_{6}, \nu_{6}) \big] = V_{4} \big( -(\lambda_{2}\nu_{6} + \nu_{2}\lambda_{6}), \ \mu_{2}\mu_{6}, \  -(\lambda_{2}\lambda_{6} + \nu_{2}\nu_{6}) \big), \\
& \big[ V_{4}(\lambda_{4}, \mu_{4}, \nu_{4}), V_{6}(\lambda_{6}, \mu_{6}, \nu_{6}) \big] = V_{2} \big( \lambda_{4}\nu_{6} + \nu_{4}\lambda_{6}, \ -\mu_{4}\mu_{6}, \ \lambda_{4}\lambda_{6} + \nu_{4}\nu_{6}) \big), \\
& \big[ V_{3}(\lambda_{3}, \mu_{3}, \nu_{3}), V_{4}(\lambda_{4}, \mu_{4}, \nu_{4}) \big] = V_{7} \big( -(\lambda_{3}\mu_{4} + \nu_{3}\lambda_{4}), \ -(\lambda_{3}\lambda_{4} + \nu_{3}\mu_{4}), \ \mu_{3}\nu_{4} \big), \\
& \big[ V_{3}(\lambda_{3}, \mu_{3}, \nu_{3}), V_{7}(\lambda_{7}, \mu_{7}, \nu_{7}) \big] = V_{4} \big( \lambda_{3}\mu_{7} + \nu_{3}\lambda_{7}, \ \lambda_{3}\lambda_{7} + \nu_{3}\mu_{7}, \ -\mu_{3}\nu_{7} \big), \\
& \big[ V_{4}(\lambda_{4}, \mu_{4}, \nu_{4}), V_{7}(\lambda_{7}, \mu_{7}, \nu_{7}) \big] = V_{3} \big( -(\lambda_{4}\mu_{7} + \mu_{4}\lambda_{7}), \ \nu_{4}\nu_{7}, \ -(\lambda_{4}\lambda_{7} + \mu_{4}\mu_{7}) \big). \\
\end{split}
\]

\end{lemm}

Set $\zeta_{i} = V_{i}(1,1,1)$ for any $1 \leq i \leq 7$.
By Lemma \ref{bracket}, we obtain Corollary \ref{bracket2} immediately.

\begin{coro} \label{bracket2}
For any $\lambda_{i}, \mu_{i}, \nu_{i} \in \mathbb{R}\ (\lambda_{i} + \mu_{i} + \nu_{i} = 0, 1 \leq i \leq 7)$,
\[
\begin{array}{llllllllll} \vspace{1mm}
\big[ \zeta_{1}, V_{5}(\lambda_{5}, \mu_{5}, \nu_{5}) \big] = -\lambda_{5} \zeta_{4}, & \big[ V_{1}(\lambda_{1},\mu_{1},\nu_{1}), \zeta_{5} \big] = -\mu_{1} \zeta_{4}, \\ \vspace{1mm}
\big[ \zeta_{2}, V_{6}(\lambda_{6}, \mu_{6}, \nu_{6}) \big] = \mu_{6} \zeta_{4}, & \big[ V_{2}(\lambda_{2},\mu_{2},\nu_{2}), \zeta_{6} \big] = \mu_{2} \zeta_{4}, \\
\big[ \zeta_{3}, V_{7}(\lambda_{7}, \mu_{7}, \nu_{7}) \big] = -\nu_{7} \zeta_{4}, & \big[ V_{3}(\lambda_{3},\mu_{3},\nu_{3}), \zeta_{7} \big] = -\mu_{3} \zeta_{4}. \\
\end{array}
\]
\end{coro}

Note that $V_{i}\ (1 \leq i \leq 7)$ is a maximal abelian subspace of $\frak{so}(7)$.
It follows that $\frak{so}(7) = V_{1} + \cdots + V_{7}$.
Let $\frak{g}_{2}$ be the Lie algebra of $G_{2}$.
Then,
\[
\begin{split}
\frak{g}_{2} &= \{ X \in \frak{so}(7) \ ;\ X(xy) = X(x)y + xX(y) \ (x,y \in \mathbb{O}) \}
\end{split}
\]
and $\frak{g}_{2}$ is spanned by $\{ V_{i}(\lambda, \mu, \nu) \ ;\ 1 \leq i \leq 7, \lambda, \mu, \nu \in \mathbb{R}, \lambda + \mu + \nu = 0 \}$.
Thus, $\dim G_{2} = 14$.
It is known that $\mathrm{rank}\, G_{2} = 2$ and $V_{i} \cap \frak{g}_{2}\ (1 \leq i \leq 7)$ is a maximal abelian subspace of $\frak{g}_{2}$.

Next, we consider the subgroups $SU(3)$ and $SO(4)$ of $G_{2}$.
For any $x \in \mathrm{Im}\mathbb{O}$, define $(G_{2})_{x} = \{ g \in G_{2} \ ;\ g(x) = x\}$.
It is known that $(G_{2})_{x} \cong SU(3)$ for any $x \in \mathrm{Im}\mathbb{O}$ with $x \not= 0$.
Hereafter, $SU(3)$ denotes $(G_{2})_{e_{1}}$.
It is also known that $G_{2}(e_{1}) = S^{6} = \{ x \in \mathrm{Im}\mathbb{O} \ ;\ |x| = 1 \}$.
Therefore, $G_{2}/SU(3) = S^{6}$.
Let $\frak{su}(3)$ be the Lie algebra of $SU(3)$.
Then, $\frak{su}(3) = (V_{1} \cap \frak{g}_{2}) + \sum_{i=2}^{7}\mathbb{R}V_{i}(0,1,-1)$.
Define $\sigma \in G_{2}$ as
\[
\sigma(e_{i}) =
\begin{cases}
 e_{i} & (i =1,2,3), \\
 -e_{j} & (j= 4,5,6,7).
\end{cases}
\]
Then, $\sigma^{2} = \mathrm{id}$ and $(G_{2})_{\sigma} = \{ g \in G_{2} \ ; \ \sigma g = g \sigma \}$ is a symmetric subgroup.
It is known that $(G_{2})_{\sigma}$ is isomorphic to $SO(4)$.
Hereafter, $SO(4)$ denotes $(G_{2})_{\sigma}$.
Let $\frak{so}(4)$ be the Lie algebra of $SO(4)$.
Then, $\frak{so}(4) = \sum_{i=1}^{3}(V_{i} \cap \frak{g}_{2})$.
Moreover, consider $SO(7)_{\sigma} = \{ g \in SO(7) \ ;\ \sigma g = g \sigma \}$.
The identity component of $SO(7)_{\sigma}$ is denoted by $SO(3) \times SO(4)$.
We now consider the $SO(4)$-action on $S^{6} = G_{2}/SU(3)$.
Then, this $SO(4)$-action is orbit equivalent to the $SO(3) \times SO(4)$-action on $S^{6}$, which is a Hermann action.
This $SO(4)$-action is hyperpolar and $\{ \cos t e_{1} + \sin t e_{4} \ ;\ t \in \mathbb{R} \}$ is a section.
Moreover, each $SO(4)$-orbit intersects $\{ \cos t e_{1} + \sin t e_{4} \ ;\ 0 \leq t \leq \pi/2 \}$ at exactly one point.

Let $B_{\frak{g}_{2}}$ denote the Killing form of $\frak{g}_{2}$.
Then, for any $X,Y \in \frak{g}_{2}$, it holds that $B_{\frak{g}_{2}}(X,Y) = 4 \mathrm{tr}(XY)$.
Set $\langle \ ,\ \rangle = (-1/8)B_{\frak{g}_{2}}$.
Then, $\langle \ ,\ \rangle$ is an invariant inner product on $\frak{g}_{2}$.
The bi-invariant Riemannian metric on $G_{2}$ induced by $\langle \ ,\ \rangle$ is denoted by the same symbol $\langle \ ,\ \rangle$.
For any $1 \leq i , j \leq 7$ and $\lambda, \lambda', \mu, \mu', \nu, \nu' \in \mathbb{R}$, where $\lambda + \mu + \nu = 0$ and $\lambda' + \mu' + \nu' = 0$,
\[
\Big\langle V_{i}(\lambda_{i}, \mu_{i}, \nu_{i}), V_{j}(\lambda', \mu', \nu') \Big\rangle = 
\begin{cases}
\lambda \lambda' + \mu \mu' + \nu \nu' & (i = j), \\
0 & (i \not= j).
\end{cases}
\]
The Ricci tensor $r_{G_{2}}$ of $G_{2}$ is given by the bi-invariant $(0,2)$-tensor induced by $-B_{\frak{g}_{2}}$ \cite{B}.
Thus, $r_{G_{2}} = 8\langle \ ,\ \rangle$.
Similarly, the Killing form $B_{\frak{so}(7)}$ of $\frak{so}(7)$ satisfies $B_{\frak{so}(7)}(X,Y) = 5 \mathrm{tr}(XY)$ for any $X,Y \in \frak{so}(7)$.
Define $\langle \ ,\ \rangle = (-1/10)B_{\frak{so}(7)}$ and the bi-invariant Riemannian metric on $SO(7)$ induced by $\langle \ ,\ \rangle$ is denoted by the same symbol.
For any $1 \leq i , j \leq 7$ and $\lambda, \lambda', \mu, \mu', \nu, \nu' \in \mathbb{R}$,
\[
\Big\langle V_{i}(\lambda_{i}, \mu_{i}, \nu_{i}), V_{j}(\lambda', \mu', \nu') \Big\rangle = 
\begin{cases}
\lambda \lambda' + \mu \mu' + \nu \nu' & (i = j), \\
0 & (i \not= j).
\end{cases}
\]
The Ricci tensor $r_{SO(7)}$ of $SO(7)$ is given by $-B_{\frak{so}(7)}$.
Thus, $r_{SO(7)} = 10 \langle \ ,\ \rangle$.


\subsection{The compact Lie group $Spin(7)$}

The results in this subsection are based on \cite{Yokota}.
For any $X_{1} \in \frak{so}(8)$, there exist unique $X_{2}, X_{3} \in \frak{so}(8)$ such that $X_{1}(a)b + aX_{2}(b) = X_{3}(ab)$ for any $a,b \in \mathbb{O}$.
We set $F_{ij} \in \frak{so}(8)$ for any $0 \leq i,j \leq 7$ as follows.
For any $x \in \mathbb{O}$,
\[
\begin{split}
& F_{i0}(x) = \frac{1}{2}(e_{i}x), \ \  F_{0i}(x) = -\frac{1}{2}(e_{i}x) \  (i \not= 0), \\
& F_{ij}(x) = \frac{1}{2}e_{j}(e_{i}x) \  (1 \leq i \not= j \leq 7).
\end{split}
\]
Then, 
\[
\begin{array}{llllllllllll}
2F_{01} = G_{01} + G_{23} + G_{45} + G_{67}, & 2F_{23} = G_{01} + G_{23} - G_{45} - G_{67}, \\
2F_{45} = G_{01} - G_{23} + G_{45} - G_{67}, & 2F_{67} = G_{01} - G_{23} - G_{45} + G_{67}, \\
2F_{02} = G_{02} - G_{13} - G_{46} + G_{57}, & 2F_{13} = -G_{02} + G_{13} - G_{46} + G_{57}, \\
2F_{46} = -G_{02} - G_{13} + G_{46} + G_{57}, & 2F_{57} = G_{02} + G_{13} + G_{46} + G_{57}, \\
2F_{03} = G_{03} + G_{12} + G_{47} + G_{56}, & 2F_{12} =  G_{03} + G_{12} - G_{47} - G_{56}, \\ 
2F_{47} = G_{03} - G_{12} + G_{47} - G_{56}, & 2F_{56} = G_{03} - G_{12} - G_{47} + G_{56}, \\
2F_{04} = G_{04} - G_{15} + G_{26} - G_{37}, & 2F_{15} = - G_{04} + G_{15} + G_{26} - G_{37}, \\ 
2F_{26} = G_{04} + G_{15} + G_{26} + G_{37}, & 2F_{37} = - G_{04} - G_{15} + G_{26} + G_{37}, \\ 
2F_{05} = G_{05} + G_{14} - G_{27} - G_{36}, & 2F_{14} = G_{05} + G_{14} + G_{27} + G_{36}, \\
2F_{27} = - G_{05} + G_{14} + G_{27} - G_{36}, & 2F_{36} = - G_{05} + G_{14} - G_{27} + G_{36}, \\ 
2F_{06} = G_{06} - G_{17} - G_{24} + G_{35}, & 2F_{17} = -G_{06} + G_{17} - G_{24} + G_{35}, \\ 
2F_{24} = -G_{06} - G_{17} + G_{24} + G_{35}, & 2F_{35} = G_{06} + G_{17} + G_{24} + G_{35}, \\
2F_{07} = G_{07} + G_{16} + G_{25} + G_{34}, & 2F_{16} = G_{07} + G_{16} - G_{25} - G_{34}, \\ 
2F_{25} = G_{07} - G_{16} + G_{25} - G_{34}, & 2F_{34} = G_{07} - G_{16} - G_{25} + G_{34}.
\end{array}
\]
Since $G_{ij}\ (0 \leq i < j \leq 7)$ form a basis of $\frak{so}(8)$, it holds that $F_{ij}\ (0 \leq i < j \leq 7)$ also form a basis of $\frak{so}(8)$.
Define the automorphisms $\alpha, \beta, \gamma$ of $\frak{so}(8)$ as follows:
For any $X \in \frak{so}(8)$ and $a \in \mathbb{O}$,
\[
\alpha(X)(a) = \overline{X(\overline{a})}, \quad \beta(G_{ij}) = F_{ij} \ (0 \leq i < j \leq 7), \quad \gamma = \beta \circ \alpha.
\]
Then, $\alpha^{2} = \beta^{2} = \mathrm{Id}$.
If $X_{1}, X_{2}, X_{3} \in \frak{so}(8)$ satisfy $X_{1}(a)b + aX_{2}(b) = X_{3}(ab)$ for any $a,b \in \mathbb{O}$, then $X_{2} = \gamma(X_{1})$ and $X_{3} = \beta(X_{1})$.
Moreover, 
\[
\frak{so}(7) = \{ X \in \frak{so}(8) \ ;\ \alpha(X) = X \}, \ \ 
\frak{g}_{2} = \{ X \in \frak{so}(8) \ ;\ \beta(X) = \gamma(X) = X \}.
\]
In particular, $X \in \frak{so}(7)$ if and only if $\beta(X) = \gamma(X)$.
For any $g_{1} \in SO(8)$, there exist $g_{2}, g_{3} \in SO(8)$ such that $(g_{1}a)(g_{2}b) = g_{3}(ab)$ for any $a,b \in \mathbb{O}$.
Such $(g_{2}, g_{3})$ is either $(g_{2}, g_{3})$ or $(-g_{2}, -g_{3})$.
Then, 
\[
\{ (g_{1}, g_{2}, g_{3}) \in SO(8)^{3} \ ;\ (g_{1}x)(g_{2}y) = g_{3}(xy) \ (x,y \in \mathbb{O}) \}
\]
is a compact Lie group isomorphic to $Spin(8)$.
Hereafter, we denote this subgroup by $Spin(8)$.
Define $\pi : Spin(8) \rightarrow SO(8) \ ;\ (g_{1}, g_{2}, g_{3}) \mapsto g_{1}$ .
Then, $\pi$ is a double covering homomorphism from $Spin(8)$ onto $SO(8)$.
Let $\frak{so}'(8)$ be the Lie algebra of $Spin(8)$.
Then,
\[
\begin{split}
\frak{so}'(8) 
&= \{ (X_{1}, X_{2}, X_{3}) \in \frak{so}(8)^{3} \ ;\ (X_{1}x)y + x(X_{2}y) = X_{3}(xy) \ (x,y \in \mathbb{O}) \} \\
&= \{ (X, \gamma(X), \beta(X)) \ ;\ X \in \frak{so}(8) \}.
\end{split}
\]
The induced automorphism from $\frak{so}'(8)$ onto $\frak{so}(8)$ by $\pi$ is denoted by the same symbol $\pi$.
Then, $\pi : \frak{so}'(8) \rightarrow \frak{so}(8)\ ;\ (X_{1}, X_{2}, X_{3}) \mapsto X_{1}$ and $\pi$ is an isomorphism.
For any $(g_{1}, g_{2}, g_{3}) \in Spin(8)$, $g_{1} \in SO(7)$ if and only if $g_{2} = g_{3}$.
Hence, the subgroup $\{ (g_{1}, g_{2}, g_{3}) \in Spin(8) \ ;\ g_{1}(1) = 1 \} = \{ (g_{1}, g_{2}, g_{3}) \in Spin(8) \ ;\ g_{2} = g_{3} \}$ is isomorphic to $Spin(7)$.
Hereafter, we denote this subgroup by $Spin(7)$.
Then, $\pi|_{Spin(7)} : Spin(7) \rightarrow SO(7)$ is a double covering homomorphism.
Let $\frak{so}'(7)$ be the Lie algebra of $Spin(7)$.
Then,
\[
\begin{split}
\frak{so}'(7) 
&= \{ (X_{1}, X_{2}, X_{3}) \in \frak{so}'(8) \ ;\ X_{1}(1) = 0 \} \\
&= \{ (X, \gamma(X), \beta(X)) \ ;\ X \in \frak{so}(8), X_{1}(1) = 1 \} \\
&= \{ (X, \gamma(X), \beta(X)) \ ;\ X \in \frak{so}(8), \beta(X) = \gamma(X) \}.
\end{split}
\]
By the definition of $G_{2}$, it is true that $\{ (g_{1}, g_{2}, g_{3}) \in Spin(8) \ ;\ g_{1} = g_{2} = g_{3} \}$ is isomorphic to $G_{2}$.
This subgroup is also denoted by $G_{2}$.
Hence, $G_{2} \subset Spin(7) \subset Spin(8)$.
Define a $Spin(7)$-representation $\chi'$ on $\mathbb{O}$ as follows.
For any $(g_{1}, g_{2}, g_{3}) \in Spin(7)$ and $a \in \mathbb{O}$,
\[
\chi'( (g_{1}, g_{2}, g_{3}))(a) = g_{2}(a).
\]
Then, $\{ g \in Spin(7) \ ;\ \chi'(g)(1) = 1 \} = G_{2}$ and $\chi'(Spin(7))(1) = S^{7} = \{ a \in \mathbb{O} \ ;\ |a| = 1 \}$.
Thus, $S^{7} = Spin(7)/G_{2}$.
The $Spin(7)$-action on $S^{7}$ induced by $\chi'$ is denoted by the same symbol.
The set of all $1$-dimensional subspaces of $\mathbb{O}$ is denoted by $\mathbb{R}P^{7}$.
The $1$-dimensional subspace spanned by $x \in \mathbb{O}\ (x \not= 0)$ is written as $[x]$.
Define an $SO(7)$-action $\chi$ on $\mathbb{R}P^{7}$ as follows.
For any $g_{1} \in SO(7)$, let $(g_{1}, g_{2}, g_{3}) \in Spin(7)$.
Then, 
\[
\chi(g_{1})[x] = [ \chi'(g_{1}, g_{2}, g_{3})(x) ] = [g_{2}(x)].
\]
Since $S^{7} = Spin(7)/G_{2}$, we obtain $\mathbb{R}P^{7} = SO(7)/G_{2}$.
For any subgroup $L \subset SO(7)$, the restriction of $\chi$ to $L$ is also denoted by the same symbol $\chi$.
Similarly, for any subgroup $L' \subset Spin(7)$, the restriction of $\chi'$ to $L'$ is also denoted by the same symbol $\chi'$.

The $G_{2}$-action $\chi$ on $\mathbb{R}P^{7}$ is orbit equivalent to the natural isometric action of the symmetric subgroup $SO(1) \times SO(7) \subset SO(8)$.
In particular, this action is a Hermann action and $\{ [\cos t + \sin t e_{4}] \ ;\ t \in \mathbb{R} \}$ is a section.
Moreover, each orbit intersects $\{ [\cos t + \sin t e_{4}] \ ;\ 0 \leq t \leq \pi/2 \}$ at exactly one point.
Next, we consider the $SO(3) \times SO(4)$-action $\chi$ on $\mathbb{R}P^{7}$.
For $\sigma \in G_{2}$, the element $(\sigma, \sigma, \sigma) \in G_{2} \subset Spin(8)$ is denoted by the same symbol.
Then, $\{ g \in Spin(7) \ ;\ \sigma g = g \sigma \}$ is a symmetric subgroup isomorphic to $Spin(3) \cdot Spin(4)$, where we use the dot product introduced by Chen and Nagano \cite{Chen-Nagano1988}.
Hereafter, $Spin(3) \cdot Spin(4)$ denotes this subgroup.
Then, $\pi(Spin(3) \cdot Spin(4)) = SO(3) \times SO(4)$.
The Lie algebra of $Spin(3) \cdot Spin(4)$ is given by $\{ X \in \frak{so}'(7) \ ;\ X \sigma = \sigma X\}$ and is spanned by $\{ ( G_{ij}, F_{ij}, F_{ij} ) \ ;\ 1 \leq i < j \leq 3 \ \text{or} \ 4 \leq i < j \leq 7 \}$.
The $Spin(3) \cdot Spin(4)$-action $\chi'$ on $S^{7}$ is orbit equivalent to the natural isometric action of the symmetric subgroup $SO(4) \times SO(4) \subset SO(8)$.
Therefore, the $SO(3) \times SO(4)$-action $\chi$ on $\mathbb{R}P^{7}$ is also orbit equivalent to the natural isometric action of $SO(4) \times SO(4) \subset SO(8)$.
This action is a Hermann action and $\{ [\cos t + \sin t e_{4}] \ ;\ t \in \mathbb{R} \}$ is a section.
Moreover, each orbit intersects $\{ [\cos t + \sin t e_{4}] \ ;\ 0 \leq t \leq \pi/2 \}$ at exactly one point.

It is known that the normalizer of $SU(3)$ in $Spin(7)$ is isomorphic to $U(3)$ (\cite{Kerr}).
This normalizer is denoted by $\tilde{U}(3)$.
The Lie algebra of $\tilde{U}(3)$ is spanned by
\[
\begin{split} &
\begin{array}{lllllllllllllll} \vspace{1mm}
\big( G_{46} + G_{57}, G_{46} + G_{57}, G_{46} + G_{57} \big), & \big( G_{47} - G_{56}, G_{47} - G_{56}, G_{47} - G_{56} \big), \\ \vspace{1mm}
\big( G_{26} + G_{37}, G_{26} + G_{37}, G_{26} + G_{37} \big), & \big( G_{27} - G_{36}, G_{27} - G_{36}, G_{27} - G_{36} \big), \\
\big( G_{24} + G_{35}, G_{24} + G_{35}, G_{24} + G_{35} \big), & \big( G_{25} - G_{34}, G_{25} - G_{34}, G_{25} - G_{34} \big), \\
\end{array} \\
& \hspace{25mm} \Big( G_{23}, F_{23}, F_{23} \Big), \ \Big( G_{45}, F_{45}, F_{45} \Big), \ \Big( G_{67}, F_{67}, F_{67} \Big).
\end{split}
\]
We easily see that $(I_{8}, -I_{8}, -I_{8}) \in \tilde{U}(3)$, where $I_{8}$ is the $8 \times 8$ identity matrix and $\pi(\tilde{U}(3)) = U(1)/\mathbb{Z}_{2} \times SU(3)$.
Hereafter, $\pi(\tilde{U}(3))$ is denoted by $U(3)$.
Let $\frak{u}(3)$ be the Lie algebra of $U(3)$.
Then, $\frak{u}(3) = V_{1} + \sum_{i=2}^{7}\mathbb{R}V_{i}(0,1,-1)$.
The $\tilde{U}(3)$-action $\chi'$ on $S^{7}$ is orbit equivalent to the natural isometric action of the symmetric subgroup $SO(2) \times SO(6) \subset SO(8)$.
Therefore, the $U(3)$-action $\chi$ on $\mathbb{R}P^{7}$ is also orbit equivalent to the natural isometric action of $SO(2) \times SO(6) \subset SO(8)$.
This action is a Hermann action and $\{ [\cos t + \sin t e_{4}] \ ;\ t \in \mathbb{R} \}$ is a section.
Moreover, each orbit intersects $\{ [\cos t + \sin t e_{4}] \ ;\ 0 \leq t \leq \pi/2 \}$ at exactly one point.
Let $\tilde{G}_{2}(\mathbb{O})$ be the set of all oriented $2$-dimensional subspaces of $\mathbb{O}$.
Then, $Spin(7)$ acts on $\tilde{G}_{2}(\mathbb{O})$ as follows:
\[
Spin(7) \times \tilde{G}_{2}(\mathbb{O}) \rightarrow \tilde{G}_{2}(\mathbb{O}) \ ;\ ((g_{1}, g_{2}, g_{3}), \langle u \wedge v \rangle ) \mapsto \langle g_{2}u \wedge g_{2}v \rangle,
\]
where $\langle u \wedge v \rangle$ is an oriented $2$-dimensional subspace spanned by $u,v$ with the orientation induced by $u, v$.
It is known that this action is transitive and the isotropy group at $\langle e_{0} \wedge e_{1} \rangle$ of $Spin(7)$ is $\tilde{U}(3)$.
Hence, $Spin(7)/\tilde{U}(3) = \tilde{G}_{2}(\mathbb{O})$ (\cite{Hashimoto-Ohashi},\cite{Kerr}).
Moreover, we see that $SO(7)$ acts on $\tilde{G}_{2}(\mathbb{O})$ via the $Spin(7)$-action since $(I_{8}, -I_{8}, -I_{8}) \langle e_{0} \wedge e_{1} \rangle  = \langle e_{0} \wedge e_{1} \rangle$.
This $SO(7)$-action is transitive and $SO(7)/U(3) = \tilde{G}_{2}(\mathbb{O})$.


\section{Type (II)}

In this subsection, we consider type (II), namely $(G,H,K) = (G_{2}, SO(4), SU(3))$.
Let $g(t) = \mathrm{exp}tV_{4}(1,-1,0) \in G_{2}\ (t \in \mathbb{R})$.
Then, $\{ g(t)(e_{1}) = \cos t e_{1} + \sin t e_{5} \ ;\ t \in \mathbb{R} \}$ is a section of the $SO(4)$-action on $S^{6} = G_{2}/SU(3)$.
Moreover, each $SO(4)$-orbit intersects $\{ \cos t e_{1} + \sin t e_{5} \ ;\ 0 \leq t \leq \pi/2 \}$ at exactly one point in $S^{6}$.
Hence, each $SO(4) \times SU(3)$-orbit intersects $S = \{ g(t) \ ;\ 0 \leq t \leq \pi/2 \}$ at exactly one point in $G_{2}$, i.e.,
\[
G_{2} = \bigcup_{0 \leq t \leq \frac{\pi}{2}}SO(4) g(t) SU(3).
\]
Note that $\{ g(t) \ ;\ t \in \mathbb{R} \}$ is not a section.
Since $T_{g(t)}Hg(t)K = \mathrm{Ad}(g(t))^{-1}\frak{h} + \frak{k}$, by using Lemma \ref{bracket}, we can verify that the tangent space $dL_{g(t)}^{-1}(T_{g(t)}Hg(t)K)$ is spanned by $\frak{su}(3)$ and 
\[
\begin{array}{lllllllllllllll} \vspace{1mm}
\cos t V_{1}(\lambda, \mu, \nu) + \sin t V_{5}(\mu, \nu, \lambda), & \cos t V_{2}(\lambda, \mu, \nu) + \sin t V_{6}( \nu, \mu, \lambda ), \\
\cos 2t V_{3}(1, 0, -1) + \sin 2t V_{7}(2, -2, 0 ), & V_{3}(1,-2,1), \\
\end{array}
\]
where $\lambda, \mu, \nu \in \mathbb{R} \ (\lambda + \mu + \nu = 0)$.
Therefore, we obtain Lemma \ref{II-dim}.

\begin{lemm} \label{II-dim}
$\dim H g(0) K = 10, \ \dim H g(t) K = 13\ (0 < t < \pi/2)$, and $\dim Hg(\pi/2)K = 11$.
\end{lemm}

Hence, if $0 < t < \pi/2$, then the $H \times K$-orbit through $g(t)$ is a principal orbit.
In this subsection, for any $1 \leq i \leq 7$, we define $S_{i} = V_{i}(0,1,-1)$ and $T_{i} = V_{i}(2,-1,-1)$.
Set $\xi = (1/\sqrt{6})dL_{g(t)}T_{4}$.
Then, $\xi$ is a unit normal vector of $Hg(t)K$ at $g(t)$.
For each $0 < t < \pi/2$, we set elements of $\frak{h} \times \frak{k}$ as follows:
\[
\begin{split}
W_{i}^{(1)} & = \frac{1}{\sqrt{2}} (0, \ -S_{i}) \ (1 \leq i \leq 7), \\
W_{1}^{(2)} & = \frac{1}{\sqrt{6}} (0, \ -T_{1}), \\
W_{2}^{(2)} & = \frac{1}{\sqrt{6}} \left( \frac{1}{\cos t}V_{2}(2,-2,0), \ -S_{2} - 2\tan t S_{6} \right), \\
W_{3}^{(2)} & = \frac{1}{\sqrt{6}} \Big( V_{3}(2,-4,2), \ -3S_{3} \Big), \\
W_{5}^{(2)} & = \frac{1}{\sqrt{6}} \left( \frac{1}{\sin t }V_{1}(-1,2,-1), \ \cot t V_{1}(-1,2,-1) \right), \\
W_{6}^{(2)} & = \frac{1}{\sqrt{6}} \left( -\frac{2}{\sin t}S_{2}, \ -2\cot t S_{2} - S_{6} \right), \\
W_{7}^{(2)} & = \frac{1}{\sqrt{6}} \left( \frac{1}{\sin 2t}V_{3}(2,0,-2) - \cot 2t V_{3}(2,-4,2), \ 4 \cot 2t S_{3} - S_{7} \right).
\end{split}
\]
Then, for any $1 \leq i,j \leq 7\ (j \not= 4)$, 
\[
(W_{i}^{(1)})^{*}_{g(t)} = \frac{1}{\sqrt{2}}dL_{g(t)}S_{i}, \quad (W_{j}^{(2)})^{*}_{g(t)} = \frac{1}{\sqrt{6}}dL_{g(t)}T_{j}.
\]
Thus, $(W_{i}^{(1)})^{*}_{g(t)}, (W_{j}^{(2)})^{*}_{g(t)}\ (1 \leq i,j \leq 7\ (j \not= 4))$ form an orthonormal basis of $T_{g(t)}Hg(t)K$.
In the following, we denote $(W_{i}^{(a)})^{*}_{g(t)}$ by $U_{i}^{a}$ and study $\langle (\nabla_{U_{i}^{a}} U_{j}^{b}), \xi \rangle$, which is denoted by $\Gamma_{(i,a)}^{(j,b)}$.
By Lemma \ref{bracket}, it follows that $\Gamma_{(i,a)}^{(j,b)} = 0$ for any $a,b = 1,2$ if $(i,j) \not\in \{ (1,5), (2,6), (3,7), (k,k) \ ;\ 1 \leq k \leq 7, k \not= 4 \}$.
Moreover, since $V_{i} \cap \frak{g}_{2}$ is a maximal abelian subspace of $\frak{g}_{2}$ for each $1 \leq i \leq 7$, 
\[
\Gamma^{(i,1)}_{(i,1)} = 0, \quad \Gamma^{(1,2)}_{(1,1)} = \Gamma^{(1,2)}_{(1,2)}
= \Gamma^{(3,2)}_{(3,1)} = \Gamma^{(3,2)}_{(3,2)} = 0
\]
for any $1 \leq i \leq 7$.
Furthermore, since $\frak{su}(3)$ is a subalgebra of $\frak{g}_{2}$ and $T_{4} \perp \frak{su}(3)$, 
\[
\Gamma^{(5,1)}_{(1,1)} = \Gamma_{(1,2)}^{(5,1)} = \Gamma_{(5,1)}^{(5,2)} 
= \Gamma^{(2,2)}_{(2,1)} = \Gamma^{(6,1)}_{(2,1)} = \Gamma^{(6,2)}_{(6,1)}
= \Gamma_{(3,1)}^{(7,1)} = \Gamma_{(7,1)}^{(7,2)} = 0. 
\]
By using Lemma \ref{bracket},
\[
\begin{array}{llllll} \vspace{1mm}
\Gamma_{(1,1)}^{(5,2)} = -\dfrac{1}{12\sqrt{2}} \langle [ S_{1}, T_{5} ], T_{4} \rangle = \dfrac{1}{2\sqrt{2}}, \\ \vspace{1mm}
\Gamma_{(1.2)}^{(5,2)} = -\dfrac{1}{12\sqrt{6}} \langle [ T_{1}, T_{5} ], T_{4} \rangle = -\dfrac{1}{2\sqrt{6}}, \\ \vspace{1mm}
\Gamma_{(5,2)}^{(5,2)} = -\dfrac{2\cot t}{6\sqrt{6}} \langle [ T_{5}, V_{1}(-1,2,-1)], T_{4} \rangle = -\dfrac{2\cot t}{\sqrt{6}}, \\ \vspace{1mm}
\Gamma_{(2,1)}^{(6,2)} = -\dfrac{1}{4\sqrt{3}} \langle [S_{2}, \frac{1}{\sqrt{6}}T_{6} - S_{6}], T_{4} \rangle = -\dfrac{1}{2\sqrt{2}}, \\ \vspace{1mm}
\Gamma_{(2,2)}^{(2,2)} = -\dfrac{\tan t}{6} \langle [T_{2}, S_{6}], T_{4} \rangle = \dfrac{2}{\sqrt{6}}\tan t, \\ \vspace{1mm}
\Gamma_{(2,2)}^{(6,1)} = -\dfrac{1}{12\sqrt{2}} \langle [T_{2}, S_{6}], T_{4} \rangle = \dfrac{1}{2\sqrt{2}}, \\ \vspace{1mm}
\Gamma_{(2,2)}^{(6,2)} = -\dfrac{1}{12\sqrt{2}} \langle [T_{2}, T_{6} - 2S_{6}], T_{4} \rangle = 0, \\ \vspace{1mm}
\Gamma_{(6,2)}^{(6,2)} = \dfrac{\cot t}{6\sqrt{6}} \langle [T_{6}, 2S_{2}], T_{4} \rangle = -\dfrac{2}{\sqrt{6}} \cot t, \\ \vspace{1mm}
\Gamma_{(3,1)}^{(7,2)} = -\dfrac{1}{12\sqrt{3}} \langle [S_{3}, T_{7} - 2S_{7}], T_{4} \rangle = \dfrac{1}{2\sqrt{2}}, \\ \vspace{1mm}
\Gamma_{(3,2)}^{(7,1)} = -\dfrac{1}{12\sqrt{2}} \langle [T_{3}, S_{7}], T_{4} \rangle = \dfrac{1}{2\sqrt{2}}, \\ \vspace{1mm}
\Gamma_{(3,2)}^{(7,2)} = -\dfrac{1}{12\sqrt{6}} \langle [T_{3}, T_{7} - 2S_{7}], T_{4} \rangle = \dfrac{2}{\sqrt{6}}, \\
\Gamma_{(7,2)}^{(7,2)} = -\dfrac{\cot 2t}{6\sqrt{6}} \langle [T_{7}, 4S_{3}], T_{4} \rangle = -\dfrac{4 \cot 2t}{\sqrt{6}}.
\end{array}
\]
Set the subspaces $L_{i} \ (1 \leq i \leq 4)$ of $T_{g(t)}Hg(t)K$ as follows:
\[
L_{i} = \mathbb{R}U_{i}^{1} + \mathbb{R}U_{i}^{2} + \mathbb{R}U_{i+4}^{1} + \mathbb{R}U_{i+4}^{2}\ (1 \leq i \leq 3), \quad L_{4} = \mathbb{R}U_{4}^{1}.
\]
The shape operator $A^{\xi}$ satisfies $A^{\xi}(L_{i}) \subset L_{i}$ for any $1 \leq i \leq 4$.
The representative matrix of $A^{\xi}|_{L_{i}}\ (1 \leq i \leq 3)$ with respect to the basis $U_{i}^{1}, U_{i}^{2}, U_{i+4}^{1}, U_{i+4}^{2}$ is given by
\[
\footnotesize 
\begin{array}{ll} \vspace{1mm}
A^{\xi}|_{L_{1}} =
\frac{1}{2\sqrt{6}}
\left(
\begin{array}{cccccc}
0 & 0 & 0 & \sqrt{3} \\
0 & 0 & 0 & -1 \\
0 & 0 & 0 & 0 \\
\sqrt{3} & -1 & 0 & -4\cot t\\
\end{array}
\right), \\ \vspace{1mm}
A^{\xi}|_{L_{2}} =
\frac{1}{2\sqrt{6}}
\left(
\begin{array}{cccccc}
0 & 0 & 0 & -\sqrt{3} \\
0 & 4 \tan t & \sqrt{3} & 0 \\
0 & \sqrt{3} & 0 & 0 \\
-\sqrt{3} & 0 & 0 & -4\cot t\\
\end{array}
\right),
\\
A^{\xi}|_{L_{3}} =
\frac{1}{2\sqrt{6}}
\left(
\begin{array}{cccccc}
0 & 0 & 0 & \sqrt{3} \\
0 & 0 & \sqrt{3} & 4 \\
0 & \sqrt{3} & 0 & 0 \\
\sqrt{3} & 4 & 0 & -8\cot 2t\\
\end{array}
\right),
\end{array}
\]
and $A^{\xi}|_{L_{4}}$ is the zero matrix.
Hence, the mean curvature vector of $Hg(t)K$ at $g(t)$ is $(2/\sqrt{6})( -2\cot t + \tan t - 2 \cot 2t)\xi$.

\begin{lemm} \label{II-minimal}
A principal orbit $Hg(t)K\ (0 < t < \pi/2)$ is minimal if and only if $t = \tan^{-1}(\sqrt{3/2})$.

\end{lemm}

We consider the principal curvatures.
For each $1 \leq i \leq 3$, set $B_{i} = 2\sqrt{6}A^{\xi}|_{L_{i}}$.
Then,
\[
\begin{split}
\det(\lambda I_{4} - B_{1}) &= \lambda^{2} ( \lambda^{2} + 4\cot t \lambda - 4), \\
\det(\lambda I_{4} - B_{2}) &= ( \lambda^{2} + 4 \cot t \lambda - 3 )( \lambda^{2} - 4 \tan t \lambda - 3), \\
\det(\lambda I_{4} - B_{3}) &= \lambda^{4} + 4(\cot t - \tan t)\lambda^{3} - 22 \lambda^{2} -12( \cot t - \tan t) \lambda + 9 \\
&= (\lambda^{2} + 4\cot t \lambda - 3)(\lambda^{2} - 4\tan t \lambda - 3).
\end{split}
\]
Hence, the principal curvatures of $Hg(t)K\ (0 < t < \pi/2)$ are
\[
\begin{split}
0, \ \ \ \frac{1}{\sqrt{6}}\tan \frac{t}{2}, \ \ & \ -\frac{1}{\sqrt{6}} \cot \frac{t}{2}, \\ 
\frac{1}{2\sqrt{6}}(2\tan t \pm \sqrt{4\tan^{2} t + 3}), \ \ & \ \frac{1}{2\sqrt{6}}(-2 \cot t \pm \sqrt{4\cot^{2} t + 3})
\end{split}
\]
with multiplicities $3,1,1,2,2,2,2$.
Therefore, the minimal principal orbit is not austere.
By Lemma \ref{biharmonic}, we obtain Lemma \ref{II-biharmonic}.

\begin{lemm} \label{II-biharmonic}
A principal orbit $Hg(t)K\ (0 < t < \pi/2)$ is proper biharmonic if and only if 
\[
t = \tan^{-1} \left( \sqrt{\frac{5 \pm \sqrt{19}}{2}} \right).
\]
\end{lemm}

Set $h(s) = \mathrm{exp}s\, T_{4}\ (s \in \mathbb{R})$.
Then, $\left\{ h(s)(e_{1}) = \cos 2s e_{1} - \sin 2s e_{5} \ ;\ s \in \mathbb{R} \right\}$ is a section of the $SO(4)$-action on $S^{6}$ and each $SO(4)$-orbit intersects $\{ h(s)(e_{1}) \ ;\ 0 \leq t \leq \pi/4 \}$ at exactly one point in $S^{6}$.
Hence, in $G_{2}$, each $SO(4) \times SU(3)$-orbit intersects $\{ h(s) \ ;\ 0 \leq s \leq \pi/4 \}$ at exactly one point.
Moreover, we see that $\{ h(s) \ ;\ s \in \mathbb{R} \}$ is a section of the $SO(4) \times SU(3)$-action on $G_{2}$.
For any $0 \leq s \leq \pi/4$, the $SO(4) \times SU(3)$-orbit through $h(s)$ coincides with the orbit through $g(2s)$.
Summarizing the arguments of this subsection and Lemma \ref{wr2}, we obtain Theorem \ref{II-main}.

\begin{thm} \label{II-main}
The $SO(4) \times SU(3)$ action on $G_{2}$ is a cohomogeneity one hyperpolar action.
Define $h(s) = \mathrm{exp}s \, V_{4}(2,-1,-1) \ (s \in \mathbb{R})$.
Then, $\{ h(s) \ ;\ s \in \mathbb{R}\}$ is a section of this action.
Moreover, each orbit intersects $\left\{ h(s)\ ;\ 0 \leq s \leq \frac{\pi}{4} \right\}$ at exactly one point.
The following properties are satisfied.

(i)\ An orbit $SO(4) h(s) SU(3)$ is singular if and only if either $s = 0$ or $s = \pi/4$.
The dimension of the orbit through $h(0)$ is $10$ and the dimension of the orbit through $h(\pi/4)$ is $11$.
These singular orbits are not totally geodesic but weakly reflective.

(ii)\ The principal curvatures of a principal orbit $SO(4)h(s)SU(3)\ (0 < s < \pi/4)$ are given in Table \ref{II-curvature}.
A principal orbit $SO(4)h(s)SU(3)$ is minimal if and only if 
\[
s = \frac{1}{2} \tan^{-1} \left( \sqrt{\frac{3}{2}} \right).
\]
A principal orbit is proper biharmonic if and only if
\[
s = \frac{1}{2} \tan^{-1} \left( \sqrt{\frac{5 \pm \sqrt{19}}{2}} \right).
\]

\begin{table}[h]
    \centering
    \begin{tabular}{c|cc|c|c|c|} \hline
principal curvatures & multiplicity \\ \hline
$0$ & $3$ \\
$\frac{1}{\sqrt{6}}\tan s$ & $1$ \\
$-\frac{1}{\sqrt{6}} \cot s$ & $1$ \\
$\frac{1}{2\sqrt{6}}(2\tan 2s + \sqrt{4\tan^{2} 2s + 3})$ & $2$ \\
$\frac{1}{2\sqrt{6}}(2\tan 2s - \sqrt{4\tan^{2} 2s + 3})$ & $2$ \\
$\frac{1}{2\sqrt{6}}(-2 \cot 2s + \sqrt{4\cot^{2} 2s + 3})$ & $2$ \\
$\frac{1}{2\sqrt{6}}(-2 \cot 2s - \sqrt{4\cot^{2} 2s + 3})$ & $2$ \\ \hline
    \end{tabular}
    \caption{The principal curvatures of the principal orbit $SO(4)h(s)SU(3)\ (0 < s < \pi/4)$}
    \label{II-curvature}
\end{table}

\end{thm}

\begin{remark}
Denote the set of all $3$-dimensional subspaces of $\mathrm{Im}\mathbb{O}$ by $G_{3}(\mathrm{Im}\mathbb{O})$.
We call $V \in G_{3}(\mathrm{Im}\mathbb{O})$ is an associative subspace if $x(yz) = (xy)z$ for any $x,y,z \in V$.
The set of all associative subspaces is denoted by $G_{ass}(\mathrm{Im}\mathbb{O})$ and called the associative Grassmannian.
Then, $G_{2}$ acts on $G_{ass}(\mathrm{Im}\mathbb{O})$ naturally and it is well known that this action is transitive \cite{HL}.
The isotropy group at $\mathbb{R}e_{1} + \mathbb{R}e_{2} + \mathbb{R}e_{3} \in G_{ass}(\mathrm{Im}\mathbb{O})$ of $G_{2}$ is $SO(4)$ and $G_{ass}(\mathrm{Im}\mathbb{O}) = G_{2}/SO(4)$.
Moreover, $G_{ass}(\mathrm{Im}\mathbb{O})$ is a totally geodesic submanifold of $G_{3}(\mathrm{Im}\mathbb{O})$.
For each $V \in G_{ass}(\mathrm{Im}\mathbb{O})$, define $\mathrm{pr}_{V} : \mathrm{Im}\mathbb{O} \rightarrow V$ as the orthogonal projection onto $V$ with respect to the inner product $(\ ,\ )$.
Then,
\[
f : G_{ass}(\mathrm{Im}\mathbb{O}) \rightarrow \mathbb{R} \ ;\ V \mapsto |\mathrm{pr}_{V}(e_{1})|
\]
is a smooth function on $G_{ass}(\mathrm{Im}\mathbb{O})$.
We see that $f$ is invariant under $SU(3)$-action and $f(G_{ass}(\mathrm{Im}\mathbb{O})) = [0,1]$.
Each level set of $f$ is an $SU(3)$-orbit.
The level sets $f^{-1}(0)$ and $f^{-1}(1)$ are the only singular orbits and $f^{-1}(0) = \mathbb{C}P^{2}$ and $f^{-1}(1) = SU(3)/SO(3)$.
These singular orbits are totally geodesic.
The other level sets are principal orbits.

\end{remark}


\section{Type (III)}

In this section, we consider type (III), namely $(G,H,K) = (SO(7), G_{2}, G_{2})$.
Let $g(t) = \mathrm{exp}tV_{4}(1,0,1)\ (t \in \mathbb{R})$.
Then, $\{ \chi(g(t))[e_{0}] = [\cos t + \sin t e_{4}] \ :\ t \in \mathbb{R} \}$ is a section of the $G_{2}$-action $\chi$ on $\mathbb{R}P^{7} = SO(7)/G_{2}$.
Moreover, each $G_{2}$-orbit intersects $S = \{ \chi(g(t))[e_{0}] \ ;\ 0 \leq t \leq \pi/2 \}$ at exactly one point.
Hence, each $G_{2} \times G_{2}$-orbit in $SO(7)$ intersects $\{ g(t) \ ;\ 0 \leq t \leq \pi/2 \}$ at exactly one point, i.e.,
\[
SO(7) = \bigcup_{0 \leq t \leq \frac{\pi}{2}}G_{2}g(t)G_{2}.
\]
Note that $\{ g(t) \ ;\ t \in \mathbb{R} \}$ is not a section.
As in Section 3, the tangent space $dL_{g(t)}^{-1}T_{g(t)}Hg(t)K$ is spanned by $\frak{g}_{2}$ and 
\[
\begin{array}{lllllllllllllll} \vspace{1mm}
\cos t V_{1}(\lambda_{1}, \mu_{1}, \nu_{1}) - \sin t V_{5}(-\mu_{1}, \lambda_{1}, \nu_{1}), \\ \vspace{1mm}
\cos t V_{5}(\lambda_{5}, \mu_{5}, \nu_{5}) - \sin t V_{1}(-\mu_{5}, \lambda_{5}, -\nu_{5}), \\ \vspace{1mm}
\cos t V_{3}(\lambda_{3}, \mu_{3}, \nu_{3}) - \sin t V_{7}(\nu_{3}, \lambda_{3}, -\mu_{3}), \\ \vspace{1mm}
\cos t V_{7}(\lambda_{7}, \mu_{7}, \nu_{7}) - \sin t V_{3}(-\mu_{7}, \nu_{7}, -\lambda_{7}), \\ \vspace{1mm}
\cos 2t V_{2}(1,0,1) - \sin 2t V_{6}(-1,0,-1) + V_{2}(0,-2,0), \\
\cos 2t V_{6}(1,0,1) - \sin 2t V_{2}(1,0,1) + V_{6}(0,-2,0), \\
\end{array}
\]
where $\lambda_{i} + \mu_{i} + \nu_{i} = 0 \ (i = 1,3,5,7)$.
Therefore, we obtain Lemma \ref{III-dim}.

\begin{lemm} \label{III-dim}
$\dim Hg(0)K = 14$ and $\dim Hg(t)K = 20\ (0 < t \leq \pi/2)$.
In particular, $Hg(0)K = G_{2}$.
\end{lemm}

Hence, if $0 < t \leq \pi/2$, then $Hg(t)K$ is a principal orbit.
Set $\zeta = (1/\sqrt{3})\zeta_{4}$ and $\xi = dL_{g(t)}\zeta$.
Then, $\xi$ is a unit normal vector of $Hg(t)K$ at $g(t)$.
For each $0 < t \leq \pi/2$, set elements of $\frak{h} \times \frak{k}$ as follows:
\[
\begin{split}
W_{k}^{(1)} &= \left( 0, \frac{1}{\sqrt{2}}V_{1}(-1,0,1) \right), \  W_{k}^{(2)} =  \left( 0, \frac{1}{\sqrt{6}}V_{k}(-1,2,-1) \right) \  (1 \leq k \leq 4),\\
W_{5}^{(1)} &= \left( 0, \frac{1}{\sqrt{2}}V_{5}(0,-1,1) \right), \  W_{5}^{(2)} = \left( 0, \frac{1}{\sqrt{6}}V_{5}(2,-1,-1) \right), \\
W_{6}^{(1)} &= \left( 0, \frac{1}{\sqrt{2}}V_{6}(-1,0,1) \right), \  W_{6}^{(2)} = \left( 0, \frac{1}{\sqrt{6}}V_{6}(-1,2,-1) \right), \\
W_{7}^{(1)} &= \left( 0, \frac{1}{\sqrt{2}}V_{7}(-1,1,0) \right), \  W_{7}^{(2)} = \left( 0, \frac{1}{\sqrt{6}}V_{7}(-1,-1,2) \right), \\
\end{split}
\]
\[
\begin{split}
W_{1}^{(3)} &= \frac{\sqrt{3}}{4} \left( \frac{1}{\sin t}V_{5}(-2,1,1), \  -\frac{1}{3}V_{1}(1,-2,1) + \cot t V_{5}(-2,1,1) \right), \\
W_{2}^{(3)} &= \frac{\sqrt{3}}{4} \left( -\frac{1}{\sin t}V_{2}(1,-2,1), \  -\frac{2 + \cos 2t}{3\sin t}V_{2}(1,-2,1) - \frac{2\cos t}{3}V_{6}(1,-2,1) \right), \\
W_{3}^{(3)} &= \frac{\sqrt{3}}{4} \left( \frac{1}{\sin t}V_{7}(1,1,-2), \  -\frac{1}{3}V_{3}(1,-2,1) + \cot t V_{7}(1,1,-2) \right), \\
W_{5}^{(3)} &= \frac{\sqrt{3}}{4} \left( -\frac{1}{\sin t}V_{1}(1,-2,1), \  -\frac{1}{3}V_{5}(-2,1,1) - \cot t V_{1}(1,-2,1) \right), \\
W_{6}^{(3)} &= \frac{\sqrt{3}}{4} \left( -\frac{1}{\sin t}V_{6}(1,-2,1), \  \frac{2\cos t}{3}V_{2}(1,-2,1) - \frac{2 + \cos 2t}{3\sin t}V_{6}(1,-2,1) \right), \\
W_{7}^{(3)} &= \frac{\sqrt{3}}{4} \left( -\frac{1}{\sin t}V_{3}(1,-2,1), \  -\frac{1}{3}V_{7}(1,1,-2) - \cot t V_{3}(1,-2,1) \right). \\
\end{split}
\]
Then, 
\[
\begin{split}
(W_{j}^{(3)})^{*}_{g(t)} &= dL_{g(t)} \frac{1}{\sqrt{3}}\zeta_{j} \ (j = 1,3,5,7), \\
(W_{2}^{(3)})^{*}_{g(t)} &= dL_{g(t)} \frac{1}{\sqrt{3}}( \sin t \zeta_{2} - \cos t \zeta_{6}), \\
(W_{6}^{(3)})^{*}_{g(t)} &= dL_{g(t)} \frac{1}{\sqrt{3}}( \cos t \zeta_{2} + \sin t \zeta_{6}).
\end{split}
\]
Hence, $(W_{i}^{(a)})_{g(t)}^{*}, (W_{j}^{(3)})_{g(t)}^{*}\ (1 \leq i,j \leq 7\ (j \not= 4), a = 1,2)$ form an orthonormal basis of $T_{g(t)}Hg(t)K$.
As in Section 3, we denote $(W_{i}^{(a)})^{*}_{g(t)}$ by $U_{i}^{a}$ and consider $\langle (\nabla_{U_{i}^{a}} U_{j}^{b}), \xi \rangle$, which is denoted by $\Gamma_{(i,a)}^{(j,b)}$.
By Lemma \ref{bracket}, if $(i,j) \not\in \{ (1,5), (5,1), (2,6), (6,2), (3,7), (7,3), (k,k) \ ;\ 1 \leq k \leq 7, k \not= 4 \}$, then $\Gamma_{(i,a)}^{(j,b)} = 0$ for any $1 \leq a, b \leq 3$.
Moreover, $\Gamma^{(4,1)}_{(4,1)} = \Gamma^{(4,1)}_{(4,2)} = \Gamma^{(4,2)}_{(4,2)} = 0$.
Since $\frak{g}_{2} \perp \zeta_{4}$, it holds that $\Gamma_{(j,b)}^{(i,a)} = 0$ for any $a,b = 1,2$ and $(i,j) \not\in \{ (1,5), (5,1), (2,6), (6,2), (3,7), (7,3), (k,k) \ ;\ 1 \leq k \leq 7, k \not= 4 \}$.
Moreover, since each $V_{k}$ is a maximal abelian subspace of $\frak{so}(7)$ for any $1 \leq k \leq 7$, 
\[
\Gamma^{(j,1)}_{(j,3)} = \Gamma^{(j,2)}_{(j,3)} = 0
\]
for any $j = 1,3,5,7$.
By Corollary \ref{bracket2},
\[
\Gamma_{(i,3)}^{(j,1)} = \Gamma_{(j,3)}^{(i,1)} = \Gamma_{(2,3)}^{(2,1)} = \Gamma_{(6,3)}^{(6,1)} = 0 \quad 
\]
for any $(i,j) = (1,5), (2,6), (3,7)$.
Similarly,
\[
\begin{split}
& \Gamma_{(1,3)}^{(5,2)} = -\Gamma_{(5,3)}^{(1,2)} = \Gamma_{(3,3)}^{(7,2)} = -\Gamma_{(7,3)}^{(3,2)} = \dfrac{1}{\sqrt{6}}, \\
& \Gamma_{(2,3)}^{(2,2)} = \Gamma_{(6,3)}^{(6,2)} = -\dfrac{\cos t}{\sqrt{6}}, \ \ \\
- &\Gamma_{(2,3)}^{(6,2)} = \Gamma_{(6,3)}^{(2,2)} = \dfrac{\sin t}{\sqrt{6}}, \ \ \\
&\Gamma_{(i,3)}^{(i,3)} = -\dfrac{\sqrt{3}\cot t}{2} \ \ (1 \leq i \leq 7, i \not= 4).
\end{split}
\]
By using Lemma \ref{bracket},
\[
\begin{split}
\Gamma_{(1,3)}^{(5,3)} &=  -\frac{1}{12} \langle [\zeta_{1}, V_{5}(4,1,1)], \zeta \rangle = 0, \\ 
\Gamma_{(3,3)}^{(7,3)} &= -\frac{1}{12} \langle [\zeta_{3}, V_{7}(1,1,4)], \zeta \rangle = 0, \\
\Gamma_{(2,3)}^{(6,3)} &=  -\frac{1}{6} \langle [ \sin t \zeta_{2} - \cos t \zeta_{6}, (\cos t \zeta_{2} + \sin t \zeta_{6})  \\
& \hspace{30mm}+ \cos t V_{2}(1,-2,1) - \frac{2 + \cos 2t}{2\sin t}V_{6}(1,-2,1)], \zeta \rangle = 0, \\
\end{split}
\]
Set subspaces $L_{i}\ (1 \leq i \leq 4)$ of $T_{g(t)}Hg(t)K$ as follows:
\[
L_{i} = \sum_{a=1}^{3}(\mathbb{R}U_{i}^{a} + \mathbb{R}U_{i+4}^{a}), \quad L_{4} = \mathbb{R}U_{4}^{1} + \mathbb{R}U_{4}^{2}.
\]
Then, $A^{\xi}(L_{i}) \subset L_{i}$ for any $1 \leq i \leq 4$.
The representative matrix of $A^{\xi}|_{L_{i}}\ (1 \leq i \leq 3)$ with respect to $U_{i}^{1}, U_{i}^{2}, U_{i+4}^{1}, U_{i+4}^{2}$ is given by
\[
\footnotesize
\begin{split}
A^{\xi}|_{L_{1}} &= A^{\xi}|_{L_{3}} = \frac{1}{2\sqrt{6}}
\left(
\begin{array}{ccc|ccccccc}
0 & 0 & 0 & 0 & 0 & 0 \\
0 & 0 & 0 & 0 & 0 & -2 \\
0 & 0 & -3\sqrt{2} \cot t & 0 & 2 & 0 \\ \hline
0 & 0 & 0 & 0 & 0 & 0 \\
0 & 0 & 2 & 0 & 0 & 0 \\
0 & -2 & 0 & 0 & 0 & -3\sqrt{2} \cot t \\
\end{array}
\right), \\
A^{\xi}|_{L_{2}} &= \frac{1}{2\sqrt{6}}
\left(
\begin{array}{ccc|ccccccc}
0 & 0 & 0 & 0 & 0 & 0 \\
0 & 0 & -2\cos t & 0 & 0 & 2\sin t \\
0 & -2\cos t & -3\sqrt{2} \cot t & 0 & -2\sin t & 0 \\ \hline
0 & 0 & 0 & 0 & 0 & 0 \\
0 & 0 & -2\sin t & 0 & 0 & -2\cos t \\
0 & 2\sin t & 0 & 0 & -2\cos t & -3\sqrt{2}\cot t \\
\end{array}
\right), \\
\end{split}
\]
and $A^{\xi}|_{L_{4}}$ is the zero matrix.
Thus, the mean curvature vector at $g(t)$ is $(-3\sqrt{3} \cot t)\xi$.

\begin{lemm} \label{III-biharmonic}
A principal orbit $Hg(t)K\ (0 < t \leq \pi/2)$ is minimal if and only if $t = \pi/2$.

\end{lemm}

We consider the principal curvatures.
Set $B_{i} = 2\sqrt{6}A^{\xi}|_{L_{i}}$.
Then,
\[
\begin{split}
\det(\lambda I_{6} - B_{1}) &= \lambda^{2}(\lambda^{2} + 3\sqrt{2} \cot t \lambda - 4)^{2} , \\
\det(\lambda I_{6} - B_{2}) &= \lambda^{4} + 6\sqrt{2} \cot t \lambda^{3} + \big( 18\cot^{2}t - 8 \big)\lambda^{2} - 24\sqrt{2} \cot t \lambda + 16 \\
&= \lambda^{2}\big( \lambda^{2} + 3\sqrt{2} \cot t \lambda - 4 \big)^{2}. \\
\end{split}
\]
Hence, the principal curvatures of $Hg(t)K\ (0 < t \leq \pi/2)$ are 
\[
0, \quad \frac{1}{4\sqrt{6}}(-3\sqrt{2}\cot t \pm \sqrt{18\cot^{2} t + 16})
\] 
with multiplicities $8,6,6$.
We see that the minimal principal orbit is austere.
By Lemma \ref{biharmonic}, Lemma \ref{III-biharmonic} follows.

\begin{lemm} \label{III-biharmonic}
A principal orbit $Hg(t)K\ (0 < t \leq \pi/2)$ is proper biharmonic if and only if $t = \cot^{-1}(4/3)$.
\end{lemm}

Set an isometry $f$ of $SO(7)$ as
\[
f : SO(7) \rightarrow SO(7) \ ;\ x \mapsto g(\pi)x^{-1}.
\]
Since $g(\pi) \in G_{2}$, it is true that $f(Hg(\pi/2)K) \subset Hg(\pi/2)K$ and $f(g(\pi/2)) = g(\pi/2)$.
Moreover, since $V_{4}$ is an abelian subspace of $\frak{so}(7)$, we have $df(\xi) = -\xi$.
Hence, the minimal orbit is weakly reflective by Lemma \ref{wr1}.
Set $h(s) = \mathrm{exp}s\,\zeta_{4}\ (s \in \mathbb{R})$.
Then, $\{ \chi(h(s))[e_{0}] = [\cos (3/2)s + \sin (3/2)s e_{4}] \ ; \ t \in \mathbb{R}\}$ is a section of the $G_{2}$-action $\chi$ on $\mathbb{R}P^{7} = SO(7)/G_{2}$.
Moreover, each $G_{2}$-orbit intersects $\{ \chi(h(s))[e_{0}] \ ;\ 0 \leq s \leq \pi/3 \}$ at exactly one point in $\mathbb{R}P^{7}$.
Hence, each $G_{2} \times G_{2}$-orbit in $SO(7)$ intersects $\{ h(s) \ ;\ 0 \leq s \leq \pi/3 \}$ at exactly one point.
Furthermore, $\{ h(s) \ ;\ s \in \mathbb{R} \}$ is a section.
For any $0 \leq s \leq \pi/3$, the $G_{2} \times G_{2}$-orbit through $h(s)$ coincides with the orbit through $g(3s/2)$.
Summarizing the arguments in this section and Lemma \ref{wr2}, we obtain Theorem \ref{III-main}.

\begin{thm} \label{III-main}
The $G_{2} \times G_{2}$ action on $SO(7)$ is a cohomogeneity one hyperpolar action.
Set $h(s) = \mathrm{exp}s \, \zeta_{4}$.
Then, $\{ h(s) \ ;\ s \in \mathbb{R} \}$ is a section and each orbit intersects $\{ h(s) \ ;\ 0 \leq s \leq \pi/3\}$ at exactly one point.
Moreover, the following properties are satisfied.

(i)\ The orbit through $h(0)$ is the only singular orbit.
The singular orbit is the subgroup $G_{2}$.
Moreover, the singular orbit is totally geodesic and weakly reflective.

(ii)\ The principal curvatures of a principal orbit $G_{2}h(s)G_{2}\ (0 < s \leq \pi/3)$ are given in Table \ref{III-curvature}.
A principal orbit $G_{2}h(s)G_{2}$ is minimal if and only if $s = \pi/3$, and this minimal orbit is weakly reflective.
A principal orbit is proper biharmonic if and only if 
\[
s = \frac{2}{3} \cot^{-1} \frac{4}{3}.
\]

\begin{table}[htbp]
    \centering
    \begin{tabular}{c|cc|c|c|c|} \hline
principal curvatures & multiplicity \\ \hline
$0$ & $8$ \\
$\frac{1}{4\sqrt{6}} \big( -3\sqrt{2}\cot \frac{3}{2}s + \sqrt{18\cot^{2} \frac{3}{2}s + 16} \big)$ & $6$ \\
$\frac{1}{4\sqrt{6}} \big( -3\sqrt{2}\cot \frac{3}{2}s - \sqrt{18\cot^{2} \frac{3}{2}s + 16} \big)$ & $6$ \\ \hline
    \end{tabular}
    \caption{The principal curvatures of the principal orbit $G_{2}h(s)G_{2} \ (0 < s \leq \pi/3)$}
    \label{III-curvature}
\end{table}

\end{thm}

\begin{remark}
We define a smooth function $f$ on $\mathbb{R}P^{7}$ as follows.
For any $[x] \in \mathbb{R}P^{7}$, we set $\mathrm{pr}_{[x]} : \mathbb{O} \rightarrow [x]$ as the orthogonal projection.
Then,
\[
f : \mathbb{R}P^{7} \rightarrow \mathbb{R} \ ;\ [x] \mapsto |\mathrm{pr}_{[x]}(e_{0})|.
\]
We see that $f$ is invariant under the $G_{2}$-action and $f(\mathbb{R}P^{7}) = [0,1]$.
The level sets $f^{-1}(1)$ and $f^{-1}(0)$ are $\{ [e_{0}] \}$ and a totally geodesic hypersurface $\mathbb{R}P^{6}$, respectively.
Then, $f^{-1}(1)$ is the only singular orbit.
The other level sets are principal orbits.

\end{remark}


\section{Type (IV)}

In this section, we consider type (IV), namely $(G,H,K) = (SO(7), SO(3) \times SO(4), G_{2})$.
Set $g(t) = \mathrm{exp}tV_{4}(1,0,0)\ (t \in \mathbb{R})$.
Then, $\left\{ \chi(g(t))[e_{0}] = \left[ \cos \frac{t}{2} + \sin \frac{t}{2} e_{4} \right] \ ;\ t \in \mathbb{R} \right\} $ is a section of the $SO(3) \times SO(4)$-action $\chi$ on $\mathbb{R}P^{7} = SO(7)/G_{2}$.
Moreover, each $SO(3) \times SO(4)$-orbit intersects $\{ \chi(g(t))[e_{0}] \ ;\ 0 \leq t \leq \pi \}$ at exactly one point.
Hence, each $(SO(3) \times SO(4)) \times G_{2}$-orbit in $SO(7)$ intersects $S = \{ g(t) \ ;\ 0 \leq t \leq \pi \}$ at exactly one point, i.e.,
\[
SO(7) = \bigcup_{0 \leq t \leq \pi}(SO(3) \times SO(4))g(t)G_{2}.
\]
Note that $\{ g(t) \ ;\ t \in \mathbb{R} \}$ is not a section.
As in Section 3, the tangent space $dL_{g(t)}^{-1}T_{g(t)}Hg(t)K$ is spanned by $\frak{g}_{2}$ and 
\[
\begin{array}{lllllllllllllll}
V_{1}(1,0,1), & \cos t V_{1}(0,1,0) - \sin t V_{5}(-1,0,0), \\
V_{2}(0,1,0), & \cos t V_{2}(1,0,1) - \sin t V_{6}(-1,0,-1), \\
V_{3}(0,1,0), & \cos t V_{3}(1,0,1) - \sin t V_{7}(1,1,0).
\end{array}
\]
Therefore, we obtain Lemma \ref{IV-dim}.

\begin{lemm} \label{IV-dim}
$\dim Hg(0)K = \dim Hg(\pi)K = 17$ and $\dim Hg(t)K = 20\ (0 < t < \pi).$
\end{lemm}

Hence, an orbit through $g(t)\ (0 < t < \pi)$ is a principal orbit.
Set $\zeta = (1/\sqrt{3})\zeta_{4}$ and $\xi = dL_{g(t)}\zeta$.
Then, $\xi$ is a unit normal vector of $Hg(t)K$ at $g(t)$.
For each $0 < t < \pi$, we set elements of $\frak{h} \times \frak{k}$ as follows:
$W_{i}^{(a)}\ (1 \leq i \leq 7, a = 1,2)$ is defined similarly to Section 4 and
\[
\begin{split}
W_{1}^{(3)} &= \frac{\sqrt{3}}{2} \left( V_{1}(1,0,1), \frac{1}{3}V_{1}(1,-2,1) \right), \\
W_{2}^{(3)} &= \sqrt{3} \left( V_{2}(0,1,0), -\frac{1}{3}V_{2}(1,-2,1) \right), \\
W_{3}^{(3)} &= \sqrt{3} \left( V_{3}(0,1,0), -\frac{1}{3}V_{3}(1,-2,1) \right), \\
W_{5}^{(3)} &= \frac{\sqrt{3}}{\sin t} \left( V_{1}( -\frac{\cos t}{2}, 1 - \frac{\cos t}{2} ), -\frac{\cos t}{2}V_{1}(1,-2,1) - \frac{\sin t}{3} V_{5}(-2,1,1) \right), \\
W_{6}^{(3)} &= \frac{\sqrt{3}}{2 \sin t} \left( V_{2}(1, -2\cos t, 1), \cos t V_{2}(1,-2,1) + \frac{\sin t}{3}V_{6}(1,-2,1) \right), \\
W_{7}^{(3)} &= -\frac{\sqrt{3}}{2 \sin t} \left( V_{3}(1, -2\cos t, 1), \cos t V_{3}(1,-2,1) - \frac{\sin t}{3} V_{7}(1,1,-2) \right).
\end{split}
\]
Then, $(W_{j}^{(3)})^{*}_{g(t)} = \frac{1}{\sqrt{3}} dL_{g(t)} \zeta_{j}$ for any $1 \leq i \leq 7\ ( i \not= 4)$ and $(W_{i}^{(a)})_{g(t)}^{*}, (W_{j}^{(3)})_{g(t)}^{*}\ (1 \leq i,j \leq 7\ (j \not= 4), a =1,2)$ form an orthonormal basis of $T_{g(t)}Hg(t)K$.
Set $U_{i}^{a}$ and $\Gamma_{(i,a)}^{(j,b)}$ as in Section 4.
Then, $\Gamma_{(i,a)}^{(j,b)}$ satisfies similar properties to those in Section 4, that is, if $(i,j) \not\in \{ (1,5), (5,1), (2,6),(6,2),(3,7),(7,3),(k,k) \ ; \ 1 \leq k \leq 7, k \not= 4\}$, then $\Gamma_{(j,a)}^{(i,b)} = 0$ for any $1 \leq a,b \leq 3$.
Moreover, $\Gamma_{(4,1)}^{(4,1)} = \Gamma_{(4,2)}^{(4,1)} = \Gamma_{(4,2)}^{(4,2)} = 0$ and $\Gamma_{(i,a)}^{(j,b)} = 0$ for any $a,b = 1,2$ and $(i,j) \in \{ (1,5), (5,1), (2,6),(6,2), $ $(3,7),(7,3),(k,k) \ ; \ 1 \leq k \leq 7, k \not= 4\}$.
Since each $V_{i}\ (1 \leq i \leq 7)$ is a maximal abelian subspace of $\frak{so}(7)$, it holds that $\Gamma_{(i,3)}^{(i,1)} = \Gamma_{(i,3)}^{(i,2)} = \Gamma_{(j,3)}^{(j,3)} = 0$ for any $1 \leq i \leq 7\ (i \not= 4)$ and $1 \leq j \leq 3$.
By Corollary \ref{bracket2}, $\Gamma_{(i,3)}^{(j,1)} = \Gamma_{(j,3)}^{(i,1)} = 0$ for any $(i,j) = (1,5), (2,6), (3,7)$ and
\[
\begin{split}
& \Gamma_{(1,3)}^{(5,2)} = - \Gamma_{(5,3)}^{(1,2)} = -\Gamma_{(2,3)}^{(6,2)} = \Gamma_{(6,3)}^{(2,2)} = \Gamma_{(3,3)}^{(7,2)} = -\Gamma_{(7,3)}^{(3,2)} = \dfrac{1}{\sqrt{6}}, \\
& \Gamma_{(j,3)}^{(j,3)} = -\sqrt{3}\cot t \ \  (j = 5,6,7).
\end{split}
\] 
By using Lemma \ref{bracket},
\[
\begin{split}
\Gamma_{(1,3)}^{(5,3)} &= -\frac{1}{6} \langle [\zeta_{1}, V_{5}(5,-1,-1)], \zeta \rangle = \frac{\sqrt{3}}{2}, \\
\Gamma_{(2,3)}^{(6,3)} &= -\frac{1}{6} \langle [\zeta_{2}, V_{6}(2,-1,2)], \zeta \rangle = \frac{\sqrt{3}}{2}, \\
\Gamma_{(3,3)}^{(7,3)} &= -\frac{1}{6} \langle [\zeta, V_{7}(2,2,-1)], \zeta \rangle = -\frac{\sqrt{3}}{2}.
\end{split}
\]
Set subspaces $L_{i} \ (1 \leq i \leq 4)$ of $T_{g(t)}Hg(t)K$ as in Section 4.
Then, $A^{\xi}(L_{i}) \subset L_{i}$ for any $1 \leq i \leq 4$.
The representative matrix of $A^{\xi}|_{L_{i}}\ (1 \leq i \leq 3)$ with respect to $U_{i}^{1}, U_{i}^{2}, U_{i+4}^{1}, U_{i+4}^{2}$ is given as follows:
\[
\footnotesize
\begin{split}
A^{\xi}|_{L_{1}} &= \frac{1}{2\sqrt{6}}
\left(
\begin{array}{ccc|cccc}
0 & 0 & 0 & 0 & 0 & 0 \\
0 & 0 & 0 & 0 & 0 & -2 \\
0 & 0 & 0 & 0 & 2 & 3\sqrt{2} \\ \hline
0 & 0 & 0 & 0 & 0 & 0 \\
0 & 0 & 2 & 0 & 0 & 0 \\
0 & -2 & 3\sqrt{2} & 0 & 0 & -6\sqrt{2} \cot t \\
\end{array}
\right), \\
\quad
A^{\xi}|_{L_{2}} &= \frac{1}{2\sqrt{6}}
\left(
\begin{array}{ccc|cccc}
0 & 0 & 0 & 0 & 0 & 0 \\
0 & 0 & 0 & 0 & 0 & 2 \\
0 & 0 & 0 & 0 & -2 & 3\sqrt{2} \\ \hline
0 & 0 & 0 & 0 & 0 & 0 \\
0 & 0 & -2 & 0 & 0 & 0 \\
0 & 2 & 3\sqrt{2} & 0 & 0 & -6\sqrt{2} \cot t \\
\end{array}
\right), \\
A^{\xi}|_{L_{3}} &= \frac{1}{2\sqrt{6}}
\left(
\begin{array}{ccc|cccc}
0 & 0 & 0 & 0 & 0 & 0 \\
0 & 0 & 0 & 0 & 0 & -2 \\
0 & 0 & 0 & 0 & 2 & -3\sqrt{2} \\ \hline
0 & 0 & 0 & 0 & 0 & 0 \\
0 & 0 & 2 & 0 & 0 & 0 \\
0 & -2 & -3\sqrt{2} & 0 & 0 & -6\sqrt{2} \cot t \\
\end{array}
\right), \\ \quad\quad 
\end{split}
\]
and $A^{\xi}|_{L_{4}}$ is the zero matrix.
Hence, the mean curvature vector of $Hg(t)K$ at $g(t)$ is $(-3\sqrt{3} \cot t)\xi$.

\begin{lemm} \label{IV-minimal}
A principal orbit $Hg(t)K\ (0 < t < \pi)$ is minimal if and only if $t = \pi/2$.
\end{lemm}

Next, we consider the principal curvatures of each principal orbit.
Note that $A^{\xi}|_{L_{1}}, A^{\xi}|_{L_{2}},$ and $A^{\xi}|_{L_{3}}$ are conjugate to each other.
Set $B_{1} = 2\sqrt{6}A^{\xi}|_{L_{1}}$.
Then,
\[
\begin{split}
\det(\lambda I_{6} - B_{1}) 
&= \lambda^{2}(\lambda^{4} + 6\sqrt{2}\cot t \lambda^{3} - 26 \lambda^{2} -24\sqrt{2} \cot \lambda + 16) \\
&= \lambda^{2}\left( \lambda^{2} + 3\sqrt{2} \Big( \cot \frac{t}{2} \Big)\lambda - 4 \right) \left( \lambda^{2} - 3\sqrt{2} \Big( \tan \frac{t}{2} \Big)\lambda - 4 \right). \\
\end{split}
\]
Hence, the principal curvatures of the principal orbit $Hg(t)K\ (0 < t < \pi)$ are
\[
0, \ 
\frac{1}{4\sqrt{6}} \Big( -3\sqrt{2} \cot \frac{t}{2} \pm \sqrt{ 18 \cot^{2} \frac{t}{2} + 16} \Big), \ 
\frac{1}{4\sqrt{6}} \Big( 3\sqrt{2} \tan \frac{t}{2} \pm \sqrt{ 18 \tan^{2} \frac{t}{2} + 16} \Big)
\]
with the multiplicity $8,3,3,3,3$.
In particular, the minimal principal orbit is austere.
By Lemma \ref{biharmonic}, we obtain Lemma \ref{IV-biharmonic}.

\begin{lemm}\label{IV-biharmonic}
A principal orbit $Hg(t)K\ (0 < t < \pi)$ is proper biharmonic if and only if 
\[
t = \cot^{-1}\left( \frac{\sqrt{14}}{6} \right) \ \text{or} \ \pi - \cot^{-1}\left( \frac{\sqrt{14}}{6} \right).
\]
\end{lemm}

Set an isometry $f$ of $SO(7)$ such that
\[
f : SO(7) \rightarrow SO(7) \ ;\ x \mapsto g(\pi/2)\sigma g(-\pi/2) x \sigma.
\]
Since $g(\pi/2)\sigma g(-\pi/2) \in SO(3) \times SO(4)$, it follows that  $f(Hg(\pi/2)K) \subset Hg(\pi/2)K$ and $f(g(\pi/2)) = g(\pi/2)$.
Moreover, we can verify that $df(\xi) = -\xi$ and the minimal orbit is weakly reflective by Lemma \ref{wr1}.
Set $h(s) = \mathrm{exp}s\,\zeta_{4}\ (s \in \mathbb{R})$.
As in Section 4, we can verify that $\{ h(s) \ ;\ s \in \mathbb{R} \}$ is a section of the $(SO(3) \times SO(4)) \times G_{2}$-action on $SO(7)$ and each orbit intersects $\{ h(s) \ ;\ 0 \leq s \leq \pi/3 \}$ at exactly one point.
Moreover, for any $0 \leq s \leq \pi/3$, the orbit through $g(3s)$ coincides with the orbit through $h(s)$.
Summarizing the arguments of this section and Lemma \ref{wr2}, we obtain Theorem \ref{IV-main}.

\begin{thm} \label{IV-main}
The $(SO(3) \times SO(4)) \times G_{2}$-action on $SO(7)$ is a cohomogeneity one hyperpolar action.
Set $h(s) = \mathrm{exp}s \, \zeta_{4}\ (s \in \mathbb{R})$.
Then, each orbit intersects $\{ h(s) \ ;\ 0 \leq s \leq \pi/3 \}$ at exactly one point and $\{ h(s) \ ;\ s \in \mathbb{R} \}$ is a section of this action.
Moreover, the following properties are satisfied.

(i)\ An orbit $(SO(3) \times SO(4))h(s)G_{2}$ is singular if and only if either $s = 0$ or $s = \pi/3$.
The dimension of these orbits is $17$.
These orbits are not totally geodesic but weakly reflective.

(ii)\ The principal curvatures of a principal orbit through $h(s)\ (0 < s < \pi/3)$ are given in Table \ref{IV-curvature}.
A principal orbit $(SO(3) \times SO(4))h(s)G_{2}$ is minimal if and only if $s = \pi/6$ and this minimal orbit is weakly reflective.
A principal orbit is proper biharmonic if and only if 
\[
s = \frac{1}{3} \cot^{-1} \left( \frac{\sqrt{14}}{6} \right) \ \text{or} \ \ s = \frac{1}{3} \left( \pi - \cot^{-1} \left( \frac{\sqrt{14}}{6} \right) \right).
\]

\begin{table}[htbp]
    \centering
    \begin{tabular}{c|cc|c|c|c|} \hline
principal curvatures & multiplicity \\ \hline
$0$ & $8$ \\
$\frac{1}{4\sqrt{6}} \Big( -3\sqrt{2} \cot \frac{3}{2}s + \sqrt{ 18 \cot^{2} \frac{3}{2}s + 16} \Big)$ & $3$ \\
$\frac{1}{4\sqrt{6}} \Big( -3\sqrt{2} \cot \frac{3}{2}s - \sqrt{ 18 \cot^{2} \frac{3}{2}s + 16} \Big)$ & $3$ \\
$\frac{1}{4\sqrt{6}} \Big( 3\sqrt{2} \tan \frac{3}{2}s + \sqrt{ 18 \tan^{2} \frac{3}{2}s + 16} \Big)$ & $3$ \\
$\frac{1}{4\sqrt{6}} \Big( 3\sqrt{2} \tan \frac{3}{2}s - \sqrt{ 18 \tan^{2} \frac{3}{2}s + 16} \Big)$ & $3$ \\ \hline
    \end{tabular}
    \caption{The principal curvatures of the principal orbit $(SO(3) \times SO(4))h(s)G_{2}\ (0 < s < \pi/3)$}
    \label{IV-curvature}
\end{table}

\end{thm}


\section{Type (V)}

In this section, we consider type (V), namely $(G,H,K) = (SO(7), U(3), G_{2})$.
Set $g(t) = \mathrm{exp}tV_{4}(0,1,1)\ (t \in \mathbb{R})$.
Then, $\{ \chi(g(t))[e_{0}] = [ \cos t + \sin t e_{4} ] \ ;\ t \in \mathbb{R} \}$ is a section of the $U(3)$-action $\chi$ on $\mathbb{R}P^{7} = SO(7)/G_{2}$.
Moreover, each $U(3)$-orbit intersects $\{ \chi(g(t))[e_{0}] \ ;\ 0 \leq t \leq \pi/2 \}$ at exactly one point.
Hence, each $U(3) \times G_{2}$-orbit in $SO(7)$ intersects $\{ g(t) \ ;\ 0 \leq t \leq \pi/2\}$ at exactly one point, i.e.,
\[
SO(7) = \bigcup_{0 \leq t \leq \frac{\pi}{2}}U(3)g(t)G_{2}.
\]
Note that $\{ g(t) \ ;\ t \in \mathbb{R} \}$ is not a section.
As in Section 3, the tangent space $dL_{g(t)}^{-1}T_{g(t)}Hg(t)K$ is spanned by $\frak{g}_{2}$ and 
\[
\begin{array}{lllllllllllllllll}
V_{1}(0,1,0), & \cos 2t V_{1}(1,0,1) - \sin 2t V_{5}(0,1,1), \\
\cos t V_{2}(0,1,-1) - \sin t V_{6}(0,1,1), & \cos t V_{6}(0,1,-1) - \sin t V_{2}(0,-1,-1), \\
\cos t V_{3}(0,1,-1) - \sin t V_{7}(0,-1,-1), & \cos t V_{7}(0,1,-1) - \sin t V_{3}(0,-1,-1).
\end{array}
\]
Therefore, we obtain Lemma \ref{V-dim}.

\begin{lemm} \label{V-dim}
$ \dim Hg(0)K = 15, \dim Hg(t)K = 20\ (0 < t < \pi/2)$, and $\dim Hg(\pi/2)K = 19$.
\end{lemm}

Hence, an orbit through $g(t)\ (0 < t < \pi/2)$ is a principal orbit.
Set $\zeta = (1/\sqrt{3})\zeta_{4}$ and $\xi = dL_{g(t)}\zeta$.
Then, $\xi$ is a unit normal vector of $Hg(t)K$ at $g(t)$.
For each $0 < t < \pi/2$, we set elements of $\frak{h} \times \frak{k}$ as follows:
\[
\begin{split}
W_{1}^{(1)} = \left( 0, -\frac{1}{\sqrt{2}}V_{1}(-1,0,1) \right), \ & W_{1}^{(2)} = \left( 0, -\frac{1}{\sqrt{6}}V_{1}(-1,2,-1) \right), \\
W_{i}^{(1)} =  \left( 0, -\frac{1}{\sqrt{2}}V_{i}(0,-1,1) \right), \ & W_{i}^{(2)} = \left( 0, -\frac{1}{\sqrt{6}}V_{i}(2,-1,-1) \right) \quad (2 \leq i \leq 7), \\
\end{split}
\]
\[
\begin{split}
W_{1}^{(3)} &= \sqrt{3} \left( V_{1}(0,1,0), \frac{1}{3}V_{1}(-1,2,-1) \right), \\
W_{2}^{(3)} &= \left( \frac{\sqrt{3}}{2\sin t} \right) \left( V_{6}(0,1,-1), \cos t V_{6}(0,1,-1) - \frac{\sin t}{3}V_{2}(2,-1,-1) \right), \\
W_{3}^{(3)} &= \left( \frac{\sqrt{3}}{2\sin t} \right) \left( V_{7}(0,1,-1), \cos t V_{7}(0,1,-1) - \frac{\sin t}{3}V_{3}(2,-1,-1) \right), \\
W_{5}^{(3)} &= \left( -\frac{\sqrt{3}}{2\sin 2t} \right) \left( V_{1}(1, -2\cos 2t, 1), -\cos 2t V_{1}(-1,2,-1) + \frac{\sin 2t}{3} V_{5}(2,-1,-1) \right), \\
W_{6}^{(3)} &= \left( -\frac{\sqrt{3}}{2\sin t} \right) \left( V_{2}(0,1,-1), \cos t V_{2}(0,1,-1) + \frac{\sin t}{3}V_{6}(2,-1,-1) \right), \\
W_{7}^{(3)} &= \left( \frac{\sqrt{3}}{2\sin t} \right) \left( V_{3}(0,1,-1), \cos t V_{3}(0,1,-1) - \frac{\sin t}{3} V_{7}(2,-1,-1) \right).
\end{split}
\]
Then, for any $1 \leq i \leq 7 \ (i \not= 4)$, it holds that $(W_{i}^{(3)})^{*}_{g(t)} = (1/\sqrt{3}) dL_{g(t)} \zeta_{i}$.
Hence, $(W_{i}^{(a)})_{g(t)}^{*}, (W_{j}^{(3)})_{g(t)}^{*}\ (1 \leq i,j \leq 7\ (j \not= 4), a =1,2)$ form an orthonormal basis of $T_{g(t)}Hg(t)K$.
Set $U_{i}^{a}$ and $\Gamma_{(i,a)}^{(j,b)}$ as in Section 4.
Then, $\Gamma_{(i,a)}^{(j,b)}$ satisfies similar properties to those in Section 4, that is, if $(i,j) \not\in \{ (1,5), (5,1), (2,6),(6,2),(3,7),(7,3),(k,k) \ ; \ 1 \leq k \leq 7, k \not= 4\}$, then $\Gamma_{(j,a)}^{(i,b)} = 0$ for any $1 \leq a,b \leq 3$.
Moreover, $\Gamma_{(4,1)}^{(4,1)} = \Gamma_{(4,2)}^{(4,1)} = \Gamma_{(4,2)}^{(4,2)} = 0$ and $\Gamma_{(i,a)}^{(j,b)} = 0$ for any $a,b = 1,2$ and $(i,j) \in \{ (1,5), (5,1), (2,6),(6,2),(3,7),(7,3),(k,k) \ ; \ 1 \leq k \leq 7, k \not= 4\}$.
Since each $V_{k}$ is a maximal abelian subspace of $\frak{so}(7)$ for any $1 \leq k \leq 7$, it holds that $\Gamma_{(i,3)}^{(i,1)} = \Gamma_{(i,3)}^{(i,2)} = \Gamma_{(1,3)}^{(1,3)} = 0$ for any $1 \leq i \leq 7 \ (i \not= 4)$.
By Corollary \ref{bracket2}, we obtain $\Gamma_{(1,3)}^{(5,1)} = \Gamma_{(5,3)}^{(1,1)} = 0$ and
\[
\begin{split}
- & \Gamma_{(1,3)}^{(5,2)} = \Gamma_{(5,3)}^{(1,2)} = \frac{1}{\sqrt{6}},\\
& \Gamma_{(5,3)}^{(5,3)} = -\sqrt{3}\cot 2t, \\
& \Gamma_{(2,3)}^{(2,3)}  = \Gamma_{(6,3)}^{(6,3)} = \Gamma_{(3,3)}^{(3,3)} = \Gamma_{(7,3)}^{(7,3)} = -\frac{\sqrt{3}}{2}\cot t, \\
- & \Gamma_{(2,3)}^{(6,1)} = \Gamma_{(6,3)}^{(2,1)} = - \Gamma_{(3,3)}^{(7,1)} = -\Gamma_{(7,3)}^{(3,1)} = \frac{1}{2\sqrt{2}}, \\
- & \Gamma_{(2,3)}^{(6,2)} = \Gamma_{(6,3)}^{(2,2)} = \Gamma_{(3,3)}^{(7,2)} = -\Gamma_{(7,3)}^{(3,2)} = \frac{1}{2\sqrt{6}}.
\end{split}
\]
By using Lemma \ref{bracket},
\[
\begin{split}
\Gamma_{(1,3)}^{(5,3)} &= -\frac{1}{6} \langle [\zeta_{1}, V_{5}(-1,2,2)], \zeta \rangle = -\frac{\sqrt{3}}{2}, \\
\Gamma_{(2,3)}^{(6,3)} &= -\frac{1}{6} \langle [\zeta_{2}, V_{6}(-1,2,2)], \zeta \rangle = 0, \\
\Gamma_{(3,3)}^{(7,3)} &= -\frac{1}{6} \langle [\zeta_{3}, V_{7}(-1,2,2)], \zeta \rangle = 0.
\end{split}
\]
Set subspaces $L_{i} \ (1 \leq i \leq 4)$ of $T_{g(t)}Hg(t)K$ as in Section 4.
Then, $A^{\xi}(L_{i}) \subset L_{i}$ for any $1 \leq i \leq 4$.
Moreover, the representative matrix of $A^{\xi}|_{L_{i}}\ (1 \leq i \leq 3)$ with respect to $U_{i}^{1}, U_{i}^{2}, U_{i+4}^{1}, U_{i+4}^{2}$ is given as follows:
\[
\begin{split}
A^{\xi}|_{L_{1}} &= 
\scriptsize
\frac{1}{2\sqrt{6}}
\left(
\begin{array}{ccc|cccc}
0 & 0 & 0 & 0 & 0 & 0 \\
0 & 0 & 0 & 0 & 0 & 2 \\
0 & 0 & 0 & 0 & -2 & -3\sqrt{2} \\ \hline
0 & 0 & 0 & 0 & 0 & 0 \\
0 & 0 & -2 & 0 & 0 & 0 \\
0 & 2 & -3\sqrt{2} & 0 & 0 & -6\sqrt{2} \cot 2t \\
\end{array}
\right), \\
\normalsize
A^{\xi}|_{L_{2}} &= 
\scriptsize
\frac{1}{2\sqrt{6}}
\left(
\begin{array}{ccc|cccc}
0 & 0 & 0 & 0 & 0 & \sqrt{3} \\
0 & 0 & 0 & 0 & 0 & 1 \\
0 & 0 & -3\sqrt{2}\cot t & -\sqrt{3} & -1 & 0 \\ \hline
0 & 0 & -\sqrt{3} & 0 & 0 & 0 \\
0 & 0 & -1 & 0 & 0 & 0 \\
\sqrt{3} & 1 & 0 & 0 & 0 & -3\sqrt{2} \cot t \\
\end{array}
\right), \\
\normalsize
A^{\xi}|_{L_{3}} &=
\frac{1}{2\sqrt{6}}
\scriptsize
\left(
\begin{array}{ccc|cccc}
0 & 0 & 0 & 0 & 0 & -\sqrt{3} \\
0 & 0 & 0 & 0 & 0 & -1 \\
0 & 0 & -3\sqrt{2}\cot t & -\sqrt{3} & 1 & 0 \\ \hline
0 & 0 & -\sqrt{3} & 0 & 0 & 0 \\
0 & 0 & 1 & 0 & 0 & 0 \\
-\sqrt{3} & -1 & 0 & 0 & 0 & -3\sqrt{2} \cot t \\
\end{array}
\right),
\end{split}
\]
and $A^{\xi}|_{L_{4}}$ is the zero matrix.
Hence, the mean curvature vector of $Hg(t)K$ at $g(t)$ is $-\sqrt{3}(\cot 2t + 2 \cot t)\xi$.

\begin{lemm}
A principal orbit $Hg(t)K\ (0 \leq t \leq \pi/2)$ is minimal if and only if $t = \tan^{-1} \sqrt{5}$.
\end{lemm}

Next, we consider the principal curvatures of each principal orbit.
$A^{\xi}|_{L_{2}}$ and $A^{\xi}|_{L_{3}}$ are conjugate to each other.
Set $B_{1} = 2\sqrt{6}A^{\xi}|_{L_{1}}$ and $B_{2} = 2\sqrt{6}A^{\xi}|_{L_{2}}$.
Then, 
\[
\begin{split}
\det (\lambda I_{6} - B_{1}) &= \lambda^{2}(\lambda^{2} + 3\sqrt{2}\cot t \lambda - 4)(\lambda^{2} - 3\sqrt{2} \tan t \lambda - 4), \\
\det (\lambda I_{6} - B_{2}) &= \lambda^{2}(\lambda^{2} + 3\sqrt{2}\cot t \lambda - 4)^{2}.
\end{split}
\]
Hence, the principal curvatures of a principal orbit $Hg(t)K\ (0 < t < \pi/2)$ are
\[
0, \  \frac{1}{4\sqrt{6}}(-3\sqrt{2} \cot t \pm \sqrt{18\cot^{2}t + 16}), \ \frac{1}{4\sqrt{6}}(3\sqrt{2} \tan t \pm \sqrt{18\tan^{2}t + 16})
\]
with the mutliplicity $8,5,5,1,1$.
Therefore, the minimal principal orbit is not austere.
By Theorem \ref{biharmonic}, we obtain Lemma \ref{V-biharmonic}.

\begin{lemm} \label{V-biharmonic}
A principal orbit $Hg(t)K\ (0 \leq t \leq \pi/2)$ is proper biharmonic if and only if 
\[
t = \tan^{-1} \left( \frac{16 \pm \sqrt{211}}{3} \right).
\]
\end{lemm}

Set $h(s) = \mathrm{exp}s\,\zeta_{4}\ (s \in \mathbb{R})$.
As in Section 4, $\{ h(s) \ ;\ s \in \mathbb{R} \}$ is a section of the $U(3) \times G_{2}$ action on $SO(7)$ and each orbit intersects $\{ h(s) \ ;\ 0 \leq s \leq \pi/3 \}$ at exactly one point.
Moreover, for any $0 \leq s \leq \pi/3$, the $U(3) \times G_{2}$-orbit through $g(3s/2)$ coincides with the orbit through $h(s)$.
Summarizing the arguments in this section and Lemma \ref{wr2}, we obtain Theorem \ref{V-main}.

\begin{thm} \label{V-main}
The $U(3) \times G_{2}$-action on $SO(7)$ is a cohomogeneity one hyperpolar action.
Set $h(s) = \mathrm{exp}s\,\zeta_{4}\ (s \in \mathbb{R})$.
Then, $\{ h(s) \ ;\ s \in \mathbb{R}\}$ is a section of this action and each $U(3) \times G_{2}$-orbit intersects $\{ h(s) \ ;\ 0 \leq s \leq \pi/3 \}$ at exactly one point.
Moreover, the following properties are satisfied.

(i)\ An orbit $U(3)h(s)G_{2}$ is singular if and only if either $s = 0$ or $s = \pi/3$.
The dimension of the orbit through $h(0)$ is $15$ and the dimension of the orbit through $h(\pi/3)$ is $19$.
These singular orbits are not totally geodesic but weakly reflective.

(2)\ The principal curvatures of a principal orbit through $h(s)\ (0 < s < \pi/3)$ are given in Table \ref{V-curvature}.
A principal orbit $U(3)h(s)G_{2}$ is minimal if and only if $s = (2/3)\tan^{-1} \sqrt{5}$.
This minimal orbit is not austere.
A principal orbit is proper biharmonic if and only if 
\[
s = \frac{2}{3} \tan^{-1} \left( \frac{16 \pm \sqrt{211}}{3} \right).
\]

\begin{table}[htbp]
    \centering
    \begin{tabular}{c|cc|c|c|c|} \hline
principal curvatures & multiplicity \\ \hline
0 & 8 \\
$\frac{1}{4\sqrt{6}}(-3\sqrt{2} \cot \frac{3}{2}s + \sqrt{18\cot^{2}\frac{3}{2}s + 16})$ & $5$ \\
$\frac{1}{4\sqrt{6}}(-3\sqrt{2} \cot \frac{3}{2}s - \sqrt{18\cot^{2}\frac{3}{2}s + 16})$ & $5$ \\
$\frac{1}{4\sqrt{6}}(3\sqrt{2} \tan \frac{3}{2}s + \sqrt{18\tan^{2}\frac{3}{2}s + 16})$ & $1$ \\
$\frac{1}{4\sqrt{6}}(3\sqrt{2} \tan \frac{3}{2}s - \sqrt{18\tan^{2}\frac{3}{2}s + 16})$ & $1$ \\ \hline
    \end{tabular}
    \caption{The principal curvatures of the principal orbit $U(3)h(s)G_{2}\ (0 < s < \pi/3)$}
    \label{V-curvature}
\end{table}

\end{thm}

\begin{remark}
In the $G_{2}$-action on $\tilde{G}_{2}(\mathbb{O}) = SO(7)/U(3)$, there exist only two singular orbits.
One is a totally geodesic complex submanifold $\tilde{G}_{2}(\mathrm{Im}\mathbb{O})$ and the other is a real form $S^{6}$.
In particular, the $G_{2}$-action is orbit equivalent to the isometric action of a symmetric subgroup $SO(2) \times SO(6) \subset SO(8)$,  which is a Hermann action.
For each $\pm V \in \tilde{G}_{2}(\mathbb{O})$, where $\pm$ indicates the orientation, we define $\mathrm{pr}_{\pm V} : \mathbb{O} \rightarrow V$ as the orthogonal projection onto $V$.
Then,
\[
f : \tilde{G}_{2}(\mathbb{O}) \rightarrow \mathbb{R} \ ;\ \pm V \mapsto |\mathrm{pr}_{V}(e_{1})|
\]
is a smooth function on $\tilde{G}_{2}(\mathbb{O})$.
We see that $f(\tilde{G}_{2}(\mathbb{O})) = [0,1]$ and $f$ is invariant under the $G_{2}$-action since $g(e_{1}) = e_{1}$ for any $g \in G_{2}$.
Then, $f^{-1}(0) = \tilde{G}_{2}(\mathrm{Im}\mathbb{O})$ and $f^{-1}(1) = S^{6}$.
The other level sets are principal orbits.

\end{remark}


\section*{Acknowledgements}

The second author was supported by JSPS KAKENHI Grant Number JP23K12980.





\end{document}